\documentclass[12pt]{amsart}
\usepackage{amsmath,graphics, enumerate, url, amssymb}
\usepackage{pst-node}
\usepackage{pst-coil}
\usepackage{pst-plot}
\usepackage{pst-grad}
\usepackage{pstricks}
\usepackage{tikz}
\usepackage{chngcntr}
\usepackage{epsfig}
\input xy
\xyoption{all}
\newtheorem{theorem}{Theorem}
\newtheorem{lemma}{Lemma}

\theoremstyle{definition}

\newtheorem*{rem*}{Remark}

\counterwithin{table}{section}

\begin{document}
\title[Zero-sum multisets mod $p$ ]{Zero-sum multisets mod $p$ with \\ an application to surface automorphisms}
\author{
Anthony Weaver}

\begin{abstract}
    Let $V_2(F_p)$ be the two-dimensional vector space over $F_p$, the field with $p$ elements, $p$ an odd prime.  We count orbits of the general linear group $GL_2(F_p)$ on certain multisets consisting of $R \geq 3$ non-zero columns from $V_2(F_p)$.   The $R$-multisets must span $V_2(F_p)$ and be `zero-sum,' that is, the  sum (mod $p$) over the columns in the multiset is $[\begin{smallmatrix} 0 \\ 0 \end{smallmatrix}]$.   The orbit count can be interpreted as the number of topological types of fully ramified actions of the elementary abelian $p$-group of rank $2$ on compact Riemann surfaces of genus $1+ Rp(p-1)/2-p^2.$     
 \end{abstract}
 
 \maketitle

\section{Introduction}

Let $p$ be an odd prime, $F_p$ the field with $p$ elements, and $V_2(F_p)$ the two-dimensional vector space over $F_p$.  The one-dimensional subspaces are the points $\overline \infty$, $\overline 0$, $\overline 1$, $\overline 2$, $\dots$, $\overline{p-1}$ of the projective line $PG(1, F_p)$.  A standard set of representatives  in $V_2(F_p)$ is given by  the correspondence
 \begin{equation}\label{E:Zpgenerators}
 \overline \infty\leftrightarrow\begin{bmatrix} 1 \\ 0 \end{bmatrix}, \quad  \overline 0 \leftrightarrow\begin{bmatrix} 0 \\ 1 \end{bmatrix}, \quad 
  \overline i\leftrightarrow\begin{bmatrix} i \\ 1 \end{bmatrix}, \quad  i=1, \dots, p-1.
 \end{equation}  
 A multiset consisting of $R\geq 3$ non-zero columns of $V_2(F_p)$  is conveniently presented as a $2 \times R$ matrix in  block form 
  \begin{equation}\label{E:pgenmatrix}
M= \begin{bmatrix} S_\infty  \mid & S_0 \mid & S_1  \mid &  & \dots &\mid & S_{p-1} \end{bmatrix},
 \end{equation}
 where $S_i$ consists of $\mu_i \geq 0$ non-zero scalar multiples of the column $\overline i$.     
The {\it zero-sum condition} on $M$ is the requirement  that both  row  sums $\equiv 0$ (mod $p$).  The tuple 
 \begin{equation}\label{E:marking}
 {\mathcal P} =(\mu_\infty, \mu_0, \mu_1, \dots, \mu_{p-1})
 \end{equation}
 is the {\it projective marking} induced by $M$.  Since we also require that $M$ have rank $2$, at least two of the $\mu$'s must be non-zero.  
 
 The general linear group $GL_2(F_p)$ acts on zero-sum multisets by column-wise left-multiplication.  Dependence relations are preserved, hence, also, the zero-sum and rank conditions.    The subscripts on  
 $$R=\mu_\infty + \mu_0 + \dots + \mu_{p-1},$$
  may be permuted, but the underlying $R$-partition remains unchanged.     
Thus  the {\it partition type} of a $GL_2(F_p)$-orbit is well-defined.    
`Forgetting' the subscripts, we may write the parts of the partition, with multiplicities,  in decreasing order of size,
\begin{equation}\label{E:abstractPartition}
{\mathcal P} = (a^{[n]},b^{[m]},\dots, z^{[t]}), \quad a>b>\dots >z  \geq 1, 
\end{equation}
and represent ${\mathcal P}$ by its {\it Ferrers diagram}, in which the parts are exhibited as rows of dots.   
The number of rows,    
$s=n+m + \dots + t,$
must satisfy 
$2 \leq s \leq p+1.$ 
(By the rank condition on $M$, $s< 2$ is inadmissible; $s>p+1$ is inadmissible because each part of the partition is associated with a distinct subspace, of which there are exactly $p+1$.)  Additionally,  if $s=2$, $\mathcal P = (R-1, 1)$ is inadmissible since the zero-sum condition cannot be satisfied. 

Invariance of the partition type under the $GL_2(F_p)$-action reduces the problem of counting orbits on zero-sum multisets  to the sub-problem of counting orbits on the collection ${\mathcal M}({\mathcal P})$ of zero-sum multisets of  a given partition type ${\mathcal P}$.
                     For $\sigma \in GL_2(F_p)$ and $M \in {\mathcal M}({\mathcal P})$, let   
     \begin{equation}\label{E:fixSet}
              \text{Fix}_{\mathcal P}(\sigma) = \{ M \in{\mathcal M}({\mathcal P}) \mid \sigma M = M \}.
              \end{equation}
              By Burnside's Lemma  (see, e.g., \cite{vLW}, Theorem~10.5),   the number of orbits of $GL_2(F_p)$ on ${\mathcal M}({\mathcal P})$ is 
             \begin{equation}\label{E:Burnside}
             \dfrac{1}{|GL_2(F_p)|} \sum_{\sigma \in GL_2(F_p)} |\text{Fix}_{\mathcal P}(\sigma)|.
             \end{equation}
   The summation can be reduced to a sum over a set of representatives of the conjugacy classes in $GL_2(F_p)$ since conjugate elements fix the same number of points.   If $\sigma$ is the identity element, 
    $|\text{Fix}_{\mathcal P}(\sigma)| = |{\mathcal M}({\mathcal P})|.$  Table~\ref{Ta:fullcount} gives  $ |{\mathcal M}({\mathcal P})|$ for each  admissible partition of $3 \leq R \leq 6.$  
 For non-identity $\sigma$,   $|\text{Fix}_{\mathcal P}(\sigma)|>0$  only if the form of ${\mathcal P}$ satisfies certain constraints involving parameters associated to the conjugacy class of $\sigma$ (Lemma~\ref{L:NEW_feas}).   A parameter of particular importance is  the {\it central indicator}, which is  the smallest power of $\sigma$ contained in the center of $GL_2(F_p)$.   
   Lemmas in Section~\ref{S:elements} enumerate the conjugacy classes in $GL_2(F_p)$ according to central indicator and eigenvalue orders (if any).  Table~\ref{Ta:elements} gives this information explicitly in the cases $p=3, 5, 7$.    
%
   If $\sigma M = M$, $M$ is a concatenation of $\sigma$-orbits which can be exhibited as rectangular arrays  overlaying or  `tiling' a Ferrers diagram of ${\mathcal P}$.   
      Enumeration of feasible orbit-tilings yields formulae for $|\text{Fix}_{\mathcal P}(\sigma)|$.  See Section~\ref{S:action}, Theorems~\ref{T:NEW_noOneEigen} and ~\ref{T:NEW_OneEigen}.
  
  In Section~\ref{S:surfaces} we describe the surface theory which motivated, and is equivalent to,  the combinatorial problem treated in the rest of the paper.    The theorem in Section~\ref{S:multiPart} is a preliminary number-theoretic result, possibly of independent interest.    In Section~\ref{S:Burnside}, we assemble the results developed in the previous sections and apply Burnside's Lemma.  In Section~\ref{S:computations} we give the complete enumeration of $GL_2(F_p)$-orbits on zero-sum $R$-multisets  in the smallest cases $p=3,5,7$, $R=3,4,5,6$.

 \section{Surface automorphisms}\label{S:surfaces}  We briefly describe the relevance of our enumeration problem to the theory of surface automorphisms.  Let ${\mathbb Z}_p^2$ be the elementary abelian $p$-group of rank $2$.  As an additive group it is isomorphic to $V_2(F_p)$ with vector addition mod $p$.   Let the group act by automorphisms on a closed surface $X$ of genus $g>1$.  In a {\it fully ramified\,} action,  the quotient surface $X/{\mathbb Z}_p^2$ has genus $0$, and the quotient map $X \rightarrow X/{\mathbb Z}_p^2$ branches over $R\geq 3$ points.   $R$ is related to $g$ by the Riemann-Hurwitz relation
 \begin{equation}\label{E:RiemHur}
g=1+ \dfrac{Rp(p-1)}{2}-p^2.
\end{equation}
    The existence of such an action is equivalent to the existence of a short exact sequence of group homomorphisms
\begin{equation}\label{E:shortexact}\notag
1 \rightarrow  \Lambda_g \overset{\iota}{\longrightarrow} \Gamma \overset{\eta}{\longrightarrow} {\mathbb Z}_p^2 \rightarrow 1,\end{equation}
where $1$ is the trivial group, and $\Gamma$ is a Fuchsian group with {\it signature} $(0;R)$ and  presentation
\begin{equation}\label{E:Gammapres}
\Gamma = \Gamma(0;R) = \langle \gamma_1, \dots \gamma_R \mid \gamma_1^p= \dots = \gamma_R^p = \prod_{i=1}^R \gamma_i = \text{id} \rangle.
\end{equation}
$\Lambda_g$ is a torsion-free ({\it surface}) group with signature $(g; - )$ and presentation
$$\Lambda_g = \langle \alpha_1, \beta_1, \dots, \alpha_g, \beta_g \mid \prod_{i=1}^g \alpha_i\beta_i\alpha_i^{-1}\beta_i^{-1} = \text {id} \rangle.$$
  Two short exact sequences determine {\it topologically equivalent\,} actions if there are group automorphisms $\phi: \Gamma \rightarrow \Gamma$ and $\psi: {\mathbb Z}_p^2 \rightarrow {\mathbb Z}_p^2$ such that the diagram
   \begin{equation}\notag
\begin{matrix}
 1&\rightarrow& \Lambda_g &\overset{i_1}\longrightarrow & \Gamma& \overset{\eta_1}{\longrightarrow}& {\mathbb Z}_p^2 &\longrightarrow &1 \\
&&\Vert&& {\bf\ \phi \downarrow} && {\bf\  \psi \downarrow}  && \\
 1 &\rightarrow& \Lambda_g &\overset{i_2}\longrightarrow & \Gamma& \overset{\eta_2}{\longrightarrow}& {\mathbb Z}_p^2 &\longrightarrow &1
\end{matrix}
\end{equation}
commutes.  
Each epimorphism $\eta_i$ specifies a {\it generating vector} 
$$(\eta_i(\gamma_1), \dots, \eta_i(\gamma_R)) \in \underbrace{{\mathbb Z}_p^2 \times \dots \times {\mathbb Z}_p^2}_{R \text{\ factors}},$$
which is an ordered $R$-tuple of elements generating the group,  whose product is the identity (this is a consequence of the final relation in \eqref{E:Gammapres}).   In an abelian group, of course, the `product'  of a set of elements is independent of the order in which they are taken, so a generating vector is just a generating multiset  which sums to the identity.     A pair
$$(\phi, \psi) \in \text{Aut}(\Gamma) \times \text{Aut}({\mathbb Z}_p^2)$$
acts on generating vectors by sending 
$$(\eta(\gamma_1), \dots, \eta(\gamma_R)) \quad \mapsto (\psi\circ\eta\circ \phi^{-1}(\gamma_1), \dots, \psi\circ\eta\circ \phi^{-1}(\gamma_R)).$$
 $\text{Aut}(\Gamma)$ merely permutes the generators $\gamma_1, \dots, \gamma_R \in \Gamma$ in all possible ways (see \cite{B} Prop. 4.2 and Remark 4.2; or \cite{BW}, Sec. 2.3).  Hence the effective action on (unordered) generating sets is by  $\text{Aut}({\mathbb Z}_p^2) =GL_2(F_p)$ alone.  The topological types of the ${\mathbb Z}_p^2$-action are thus in bijection with $GL_2(F_p)$-orbits on zero-sum $R$-multisets of columns in $V_2(F_p)$.  Table~\ref{Ta:topTypes} summarizes the results of Section~\ref{S:computations} from this point of view.
 
The notion of topological equivalence of group actions goes back at least to Nielsen  \cite{N27}, who studied the problem for cyclic groups of prime order.   This early work was extended by Harvey \cite{H} and others in the early 1970's.   Since then there have been many attempts to count topological types of finite group actions in a given genus.  The main motivation is the bijection between topological types and conjugacy classes of finite subgroups of the mapping class group.  Lloyd \cite{L} gave a generating function for the number of topologically inequivalent ${\mathbb Z}_p$-actions in genus $g$.   Gilman \cite{G76} gave an explicit count in this cyclic case, using combinatorial techniques similar to ours.  Costa and Natanzon \cite{CN} gave a bijective classification of topological types of elementary abelian $p$-group actions which specifies (without solving) the crucial counting problem.  Broughton and Wootton \cite{BW}  describe a technique using representation theory and M\"obius inversion to count topological types of abelian actions, and provide explicit results in some low genus cases, e.g., actions of ${\mathbb Z}_p^2$ with four branch points (cf. their Section~4).  Their method depends on enumeration of the finite subgroups of the symmetric groups.  

 \section{Weighted multi-partitions}\label{S:multiPart}  
 The theorem in this section is used several times in subsequent sections.  The proof, which is somewhat involved, is given in the Appendix.   
  
 Let $(\mu_1, \mu_2, \dots, \mu_m)$ be an $m$-tuple of positive integers,  and ($\omega_1, \omega_2, \dots, \omega_m)$ an $m$-tuple of {\it weights}, which are non-zero residue classes mod $p$.  Let  $\alpha$ be an arbitrary residue  class mod $p$.    We seek $(\mu_1 + \dots + \mu_m)$-tuples
\begin{equation}\label{E:tuple}
(r_1, \dots, r_{\mu_1}, s_1, \dots, s_{\mu_2}, \dots, z_1, \dots, z_{\mu_m})
\end{equation} 
of non-zero residue classes mod $p$ which satisfy the congruence
  \begin{equation}\label{E:pseudoPart}
\omega_1(r_1  + \dots + r_{\mu_1}) + \omega_2(s_1  + \dots + s_{\mu_2}) + \dots + \omega_m(z_1 + \dots + z_{\mu_m}) \equiv \alpha \pmod p.
\end{equation}
The sums in parentheses are to be interpreted as unordered sums, or {\it partitions}, so we may assume    
\begin{gather}\label{E:partCond}
p-1 \geq r_1 \geq \dots \geq r_{\mu_1} \geq 1, \\ \notag
 p-1 \geq s_1  \geq \dots \geq s_{\mu_2} \geq 1,  \\ \notag
\vdots \\ \notag
p-1 \geq z_1  \geq \dots \geq z_{\mu_m} \geq 1.
 \end{gather}
 We call \eqref{E:pseudoPart}, subject to \eqref{E:partCond},  a  {\it weighted multi-partition congruence} and a solution \eqref{E:tuple} a {\it weighted multi-partition}. The number of  solutions   is independent of the weights.  Indeed \eqref{E:tuple} is a solution to \eqref{E:pseudoPart} if and only if
$$(\omega_1^{-1}r_1, \dots, \omega_1^{-1}r_{\mu_1}, \omega_2^{-1}s_1, \dots, \omega_2^{-1}s_{\mu_2}, \dots, \omega_m^{-1}z_1, \dots, \omega_m^{-1}z_{\mu_m})$$
is a solution to 
$$(r_1  + \dots + r_{\mu_1}) + (s_1  + \dots + s_{\mu_2}) + \dots +(z_1 + \dots + z_{\mu_m}) \equiv \alpha \pmod p,$$
with all weights equal to $1$.  Inverses are taken mod $p$,    and sub-tuples such as $\omega_1^{-1}r_1, \dots, \omega_1^{-1}r_{\mu_1}$ are reordered, if necessary,  to satisfy \eqref{E:partCond}.   
Perhaps surprisingly,  the number of solutions of \eqref{E:pseudoPart} is not always independent of the residue class  $\alpha$.  The theorem below shows that solutions are equidistributed across the residue classes if and only if the $m$-tuple $(\mu_1, \dots, \mu_m)$ contains an integer $\not\equiv 0, 1$ mod $p$.

\begin{theorem}\label{T:multiPart} Let $(\mu_1, \mu_2, \dots, \mu_m)={\mathcal P}$ be an $m$-tuple of positive integers, ($\omega_1, \omega_2, \dots, \omega_m)$ an $m$-tuple of non-zero residue classes mod $p$, and $\alpha$  an arbitrary residue  class mod $p$.  
  Let $n_{{\mathcal P}, \alpha}$ be the number of $(\mu_1 + \dots + \mu_m)$-tuples \eqref{E:tuple}  satisfying  the weighted multi-partition congruence \eqref{E:pseudoPart}, subject to \eqref{E:partCond}.   Let 
  $$B_{\mathcal P} =\prod_{i=1}^m  \binom{\mu_i+p-2}{\mu_i}.$$
  Then   
    \begin{equation}\notag
(n_{{\mathcal P}, 0}, \ n_{{\mathcal P}, 1},\   \dots,\  n_{{\mathcal P}, p-1}) = 
  (W_{\mathcal P}, Z_{\mathcal P}, \dots, Z_{\mathcal P})  \end{equation} 
 where 
 \begin{align}\notag
 W_{\mathcal P} &=
 \begin{cases}
 B_{\mathcal P}/p & \text{if $\exists i, \mu_i \not\equiv 0, 1 \pmod p$}; \\
 \dfrac{B_{\mathcal P} + (-1)^{t+1}}{p} + (-1)^t & \text{otherwise;}
 \end{cases}
 \\ \notag
 Z_{\mathcal P} &=
 \begin{cases}
 B_{\mathcal P}/p & \qquad\qquad\text{if $\exists i, \mu_i \not\equiv 0, 1 \pmod p$}; \\
 \dfrac{B_{\mathcal P} + (-1)^{t+1}}{p}  &\qquad\qquad \text{otherwise;}
 \end{cases}
 \end{align}
 and  $t\geq 0$ is the number of elements $\mu_i 
\in {\mathcal P}$ such that  $\mu_i \equiv 1$ mod $p$.
\end{theorem}
If the $m$-tuple  is a singleton $(\mu_1)$,  the multi-partition \eqref{E:pseudoPart} reduces to an ordinary partition.  On the other hand, if the $m$-tuple consists of $m$  $1$'s, \eqref{E:pseudoPart} is an ordered additive {\it composition} with $m$ positive summands.  For example, if $p=7$ and $(\mu_1)=(3)$, we count $8$ ordinary $3$-part partitions of $\alpha \equiv 0$ mod $p$, namely
\begin{equation}\label{E:7Ex}
 \begin{array}{ccccccc}
 1+2+4 && 1+1+5 && 2+2+3 && 3+3+ 1 \\
 3+5+6 && 4+4 +6 && 5+5+4 && 6 + 6 + 2,
 \end{array}
 \end{equation}
 while if $(\mu_1, \mu_2, \mu_3)=(1,1,1)$,  we count $30$ $3$-part compositions, since each distinguishable reordering of the  eight partitions contributes to the count.  For the intermediate $3$-part case, 
 $$(r_1 + r_2) + s \equiv 0 \pmod 7,$$
  where $(\mu_1, \mu_2) = (2,1)$, we assume $r_1 \geq r_2$, so the six sums in \eqref{E:7Ex} of the form $a+a+b$ count twice as $(a+a) + b$ and $(a+b) + a$, while the two sums of the form $a+b+c$ count three times each as $(a+b) + c$, $(a+c) + b$, and $(b+c) + a$, for a total of $18$.

  It is convenient, in the proof of the theorem, and elsewhere in the paper, to interpret the binomial coefficient 
\begin{equation}\label{E:ballsInBoxes}
\binom{N+ M-1}{N}
\end{equation}
  as the number of ways of distributing $N$ identical balls in $M$ distinguishable boxes.   The contents of the individual boxes are separated by the insertion of $M-1$ barriers, arbitrarily,  in a linear array of $N$ balls.

\section{The cardinality of ${\mathcal M}({\mathcal P})$}\label{S:rawCount}  
To determine $|{\mathcal M}({\mathcal P})|$ for a general partition ${\mathcal P}$, we begin with  types
\begin{equation}\label{E:allones}
{\mathcal P} = (1^{[s]}) =  
(\underbrace{1,  \dots, 1}_{s}),
\quad 2 \leq s  \leq p+1.
\end{equation}   
Let $A$ be a $2 \times s$ matrix consisting of $s$ distinct columns chosen from 
the set of standard representatives of $PG(1,F_p)$,
\begin{equation}\label{E:projLine}
\biggl\{ \begin{bmatrix} 1 \\ 0 \end{bmatrix}, \begin{bmatrix} 0 \\ 1 \end{bmatrix}, 
\begin{bmatrix} 1 \\ 1 \end{bmatrix}, \begin{bmatrix} 2 \\ 1 \end{bmatrix},\dots, \begin{bmatrix} p-1 \\ 1 \end{bmatrix} \biggr\}.
\end{equation}A zero-sum multiset  $M=M(A) \in {\mathcal M}({\mathcal P})$  corresponds to a solution $\overline x = (x_1, x_2, \dots, x_s)$ in positive  integers $1 \leq x_i \leq p-1$, of the homogeneous system 
\begin{equation}\label{E:Asystem}
A\overline x \equiv [\begin{smallmatrix} 0 \\ 0 \end{smallmatrix}]\pmod p.
\end{equation}
There are no such solutions if $s=2$ or $1$.  
If $s \geq 3$, let $d_s$ be the number of solutions.  We could arbitrarily choose $s-2$ of the $x_i$'s  in any of $(p-1)^{s-2}$ ways,  and obtain  solutions $(x_1, x_2, \dots, x_s)$ mod $p$ by solving for the remaining two $x_i$'s.  If $s=3$, all the  $x_i$'s obtained are  non-zero, and this yields $d_3=p-1$.     If $s>3$,  there is no guarantee that the two  coordinates found in this way are both positive (i.e., $\not\equiv 0$ mod $p$).
To get around this obstacle,  note that if $s \leq p-1$ we may assume, by a change of basis, that  $A$ does not contain the columns $[\begin{smallmatrix} 1 \\ 0 \end{smallmatrix}]$ or $[\begin{smallmatrix} 0 \\ 1 \end{smallmatrix}]$.  Then $A$ has the form
$$\begin{bmatrix} a_1 & a_2 & \dots & a_s \\
1 & 1 & \dots & 1 \end{bmatrix},\quad 1 \leq a_{1} <a_{2} < \dots<a_{s} \leq p-1.$$
If $\ast$, $\ast\ast$ denote arbitrary residue classes mod $p$, not necessarily distinct, the number of positive solutions of 
$$A\overline x \equiv [\begin{smallmatrix} \ast \\ 0 \end{smallmatrix}]\pmod p$$
is equal to the number of positive solutions of
\begin{equation}\label{E:latter}
A\overline x \equiv [\begin{smallmatrix} 0 \\ \ast\ast \end{smallmatrix}]\pmod p,
\end{equation}
since $(x_1, x_2, \dots, x_s)$ is a positive solution of the former if and only if $(a_1^{-1}x_1, a_2^{-1}x_2, \dots, a_s^{-1}x_s)$ is a positive solution of the latter.   By Theorem~\ref{T:multiPart}, the number of positive solutions of \eqref{E:latter} is 
\begin{equation}\label{E:ws} 
 W_{(1^{[s]})} = \frac{(p-1)^s + (-1)^{s+1}}{p} + (-1)^s, \quad s \geq 1.
\end{equation}

With \eqref{E:ws}  and the principle of inclusion and exclusion, 
we derive a recursion for $d_s$ as follows.  Let $A'$ consist of any $s-2$ columns of $A$, say, for definiteness, the last $s-2$ columns.  Choose positive $\overline{y}=(x_3, \dots, x_s)$ in all possible ways, obtaining a set    
of cardinality $(p-1)^{s-2}$.  We count those $\overline y$'s for which   $A' \overline y \equiv [\begin{smallmatrix}\alpha \\ \beta \end{smallmatrix}]$, $\alpha, \beta \not\equiv 0$ mod $p$.   Each such $\overline y$ can be uniquely completed to a positive solution $(x_1, x_2, x_3, \dots, x_s)$, $1 \leq x_i \leq p-1$, of the original system \eqref{E:Asystem}, by assigning $x_1 = -(a_{13}x_3 + a_{14} x_4 +\dots + a_{1s} x_s) \equiv \alpha \pmod p$ and $x_2 = -(x_3 + x_4 + \dots + x_s) \equiv \beta \pmod p$.  From the set of $(p-1)^{s-2}$ $\overline y$'s, we first discard the $W_{(1^{[s-2]})}$ $\overline y$'s making the first row $\equiv 0$ mod $p$.  We then separately discard those $\overline y$'s making the second row $\equiv 0$.  There are $W_{(1^{[s-2]})}$ of these $\overline y$'s as well.  
To the remaining set of $(p-1)^2 - 2W_{(1^{[s-2]})}$ $\overline y$'s, we restore the $d_{s-2}$ $\overline y$'s making both rows   {\it simultaneously} $\equiv 0$.  It follows that $d_s = (p-1)^{s-2} - 2W_{(1^{[s-2]})}+ d_{s-2},\quad s\geq 3$.  From this recursion we obtain the explicit formula
\begin{align}\label{E:recursion} 
d_s &= (p-1) \sum_{k=0}^{s-3}(-1)^k\binom{s-1}{k} p^{s-3-k} \\ \notag
&=\frac{p-1}{p^2}\biggl((p-1)^{s-1} + (-1)^{s-1}((s-1)p -1)\biggr),\qquad 3 \leq s \leq p+1.
\end{align}
It is convenient to extend \eqref{E:recursion} by defining
 
\begin{equation}\label{E:extRecursion}
d_s = 
\begin{cases}
1 & \text{if $s=0$} \\
0 & \text{if $s=1,2$.} 
\end{cases}
\end{equation}
The first few cases are
\begin{alignat}{2}\notag
d_0&=1; \\ \notag
d_1 &= 0; \\ \notag
d_2 &= 0; \\ \notag
d_3 &=(p-1);  \\ \notag
d_4 &=
(p-1)(p-3);  \\ \notag
d_5 &=
(p-1)(p^2-4p+6); \\ \notag
d_6 &=
(p-1)(p^3-5p^2+10p-10); \\ \notag
d_7 &= (p-1)(p^4 - 6p^3 + 15p^2 - 20 p + 15); \\ \notag
d_8 &= (p-1)(p^5-7p^4 +21p^3-35p^2 + 35p-21).
\end{alignat}

\begin{theorem}\label{T:all1} Let  
${\mathcal P} = (1^{[s]}), s \leq p+1$.  
\begin{equation}\notag
|{\mathcal M}({\mathcal P})| = \binom{p+1}{s} d_s.
\end{equation}
\end{theorem}
\begin{proof} We need only note that the $s$ columns comprising $A$ can be selected from $PG(1, F_p)$ in $\binom{p+1}{s}$ ways. 
\end{proof}

  We now consider more general partition types of the form ${\mathcal P} =(a^{[n]},b^{[m]},1^{[s]} )$, $a>b >1$, $n,m,s  \geq 0$.   This will suffice for our computations, but Theorem~\ref{T:generalP} below,  which treats this case, can be extended in a straightforward way to the general case.  The theorem reduces to Theorem~\ref{T:all1} in the case $n=m=0$.

\begin{theorem}\label{T:generalP}
Let
${\mathcal P} = (a^{[n]},b^{[m]},1^{[s]})$, $a>b >1$, $n,m,s \geq 0$, $n+m+s \leq p+1$.
  \begin{equation}\notag
 |{\mathcal M}({\mathcal P})| = \frac{(p+1)!}{s!(p+1-s-n-m)!}\sum_{i=0}^n\sum_{j=0}^m \frac{W_{(a)}^{n-i}Z_{(a)}^iW_{(b)}^{m-j}Z_{(b)}^j\ d_{s+i+j}}{i!j!(n-i)!(m-j)!},
  \end{equation}
  where $d_s$ is defined at  \eqref{E:extRecursion} and $W_{(\ast)}$, $Z_{(\ast)}$ are defined in Theorem~\ref{T:multiPart}.
     \end{theorem}
\begin{proof}
  We construct zero-sum multisets by a process of expansion  starting from solutions of \eqref{E:Asystem}, 
   where the coefficient matrix consists of  $s'$ distinct columns chosen from \eqref{E:projLine},  $s \leq s'\leq s+n+m\leq p+1$.  By Theorem~\ref{T:all1} there are $d_{s'}$ such solutions.  Associate entries from ${\mathcal P}$ with the elements of a solution $\overline x =(x_1, x_2, \dots,  x_s, \dots x_{s'})$  as follows. 
     Partition the set $\{x_1, \dots, x_{s'}\}$ into three disjoint subsets:  $X_a$, $X_b$ and $X_1$, consisting of $i$, $j$ and $s$  elements, respectively, where $i \leq n$, $j \leq m$, and $i+j = s'-s$.  Associate $i$ entries equal to $a$ with the elements of $X_a$, $j$ entries equal to $b$ with the elements of $X_b$, and all $s$ entries equal to $1$ with the elements of $X_1$.   The remaining $n-i$ entries equal to $a$ and $m-j$ entries equal to $b$, if any, are  associated to columns in \eqref{E:projLine} outside of $A$.  
  For each $a$  assigned to an $x_k$,   expand $x_k \not\equiv 0$ to any of the $Z_{(a)}$ $a$-part partitions of $x_k$ with  parts of size at most $p-1$.  Similarly, for each $b$ assigned to an $x_l$,  expand  $x_l$ by any of the $Z_{(b)}$ $b$-part partitions of $x_l \not\equiv 0$ with parts of size at most $p-1$.    The remaining $n-i$ $a$'s (resp. $m-j$ $b$'s) are expanded to any of the $W_{(a)}$ $a$-part (resp., any of the $W_{(b)}$ $b$-part partitions) of an integer $\equiv 0$ mod $p$. 
     
     We count the independent choices at each step of the foregoing process.  First choose $s+n+m$ of the $p+1$ columns of \eqref{E:projLine},
   $$ \binom{p+1}{n+m+s};$$
  choose $s'$ columns from this set,
  $$\binom{n+m+s}{s'};$$
  distribute $a$-, $b$- and $1$- parts of ${\mathcal P}$ over these $s'$ columns in
  $$\frac{s'!}{s!i!j!}$$
  ways; distribute any leftover $a$- and $b$-parts in
  $$\frac{(n+m + s - s')!}{(n-i)!(m-j)!}$$
  ways.
  The product of these numbers is 
  $$\frac{(p+1)!}{s!(p+1-s-n-m)!i!j!(n-i)!(n-j)!}.
  $$
  Each distribution of parts is associated with $d_{s+i+j}$ solutions of \eqref{E:Asystem}, each of which can be expanded in $W_{(a)}^{n-i}Z_{(a)}^iW_{(b)}^{m-j}Z_{(b)}^j$ ways.  Summing over all possible $i$ and $j$ yields the formula.  
   \end{proof}

Results of Theorem~\ref{T:generalP} are given in Table~\ref{Ta:fullcount} for all admissible partitions of $3 \leq R \leq 6$.  Formulae given for $p>3$ are valid for $p=3$ as well, except in the cases where the partition contains a part $3$ or $4$ ($\equiv 0, 1$ mod $3$).  Recall that the inadmissible partition types are the singletons $(R^{[1]})$ and the two-part partitions of the form $((R-1)^{[1]}, 1^{[1]})$.
 \begin{table}[h]
      \begin{center} 
       \begin{tabular}{| c | c || c | c |  }
       \hline
    \multicolumn{2}{| c ||}{$R$-partition} & \multicolumn{2}{c |}{$|{\mathcal M}({\mathcal P})|$} 
      \\ \hline
    $R$ & ${\mathcal P}$ & $p=3$ & $p>3$ 
      \\ \hline\hline
        $3$ & $(1^{[3]})$ & $8$ & $(p+1)p(p-1)^2/6 $    \\
        \hline
        $4$ & $(1^{[4]})$ & $0$ & $(p+1)p(p-1)^2(p-2)(p-3)/24 $   \\
        
        & $(2^{[1]},1^{[2]})$ &$24$ & $(p+1)p(p-1)^3/4 $ \\
               & $(2^{[2]})$ &$6$ & $(p+1)p(p-1)^2/8 $  \\
               
        \hline 
        $5$ & $(1^{[5]})$ & $0$ & $(p+1)p(p-1)^2(p-2)(p-3)(p^2-4p+6)/120 $  \\
             
             & $(2^{[1]},1^{[3]})$ &$8$ & $(p+1)p(p-1)^3(p-2)^2/12 $  \\
             
             & $(2^{[2]},1^{[1]})$ & $24$ & $(p+1)p(p-1)^4/8 $  \\
         
         & $(3^{[1]},1^{[2]})$ & $24$& $(p+1)^2p(p-1)^3/12 $ \\
         & $(3^{[1]},2^{[1]})$ & $24$& $(p+1)^2p(p-1)^2/12 $  \\
         
         \hline
         $6$ & $(1^{[6]})$ & $0$ & $(p+1)p(p-1)^2(p-2)(p-3)(p-4)(p^3-5p^2 +10p-10)/720$  \\
        
         & $(2^{[1]},1^{[4]}) $ & $0$ & $(p+1)p(p-1)^3(p-2)(p-3)(p^2-3p+3)/48 $  \\
         
         & $(2^{[2]},1^{[2]})$ &$24$ &$(p+1)p(p-1)^5(p-2)/16 $  \\
         
         & $(3^{[1]},1^{[3]})$ &$16$ & $ (p+1)^2p(p-1)^3(p-2)^2/36$ \\
        
         & $(4^{[1]},1^{[2]})$ & $48$  & $(p+2)(p+1)^2p(p-1)^3/48 $ \\
         & $(3^{[1]},2^{[1]},1^{[1]})$ & $48$  & $(p+1)^2p(p-1)^4/12 $ \\
         
          & $(2^{[3]})$ & $12$ & $(p+1)p^2(p-1)^4/48$ \\

         & $(3^{[2]})$ & $24$ & $ (p+1)^3p(p-1)^2/72$ \\
                  
         & $(4^{[1]},2^{[1]})$ & $12$ & $(p+2)(p+1)^2p(p-1)^2/48 $ \\
         \hline
        \end{tabular}
    \caption{$|{\mathcal M}({\mathcal P})|$ for admissible partitions ${\mathcal P}$ of $3 \leq R \leq 6$.}\label{Ta:fullcount}
    \end{center}
    \end{table}
 
   \section{Conjugacy classes in $GL_2(F_p)$}\label{S:elements}
   
$GL_2(F_p)$  has order $(p^2-1)(p^2-p)$.   Let $\sigma \in GL_2(F_p)$ be a non-central element, and let $\Tilde{\sigma}$ denote its image in the projective group  $PGL_2(F_p)=GL_2(F_p)/ \{[\begin{smallmatrix} x & 0 \\ 0 & x \end{smallmatrix}]: x \in F_p^\ast\}$.  $\sigma$ is {\it elliptic}, {\it parabolic}, {\it hyperbolic} if 
$\Tilde{\sigma}$ has, respectively, zero, one, or two fixed points on $PG(1, F_p)$.  Equivalently, $\sigma$ has, respectively, zero, one or two eigenspaces in $V_2(F_p)$.   Representatives of each type are given in  Table~\ref{Ta:classes} (see \cite{FH}, \S 5.2).  In the Table, $\epsilon$ denotes a primitive generator of the multiplicative group $F_p^\ast$ of non-zero elements in $F_p$.    Each geometric type is a union of conjugacy classes; the number of classes is the number of distinct representatives of the given type.  The factor of $\frac 12$ for elliptic elements comes from the fact that representatives in which $y$ is replaced by $-y$ are conjugate.  Similarly the factor of $\frac 12$ for hyperbolic elements comes from the fact that representatives with $x$ and $y$ interchanged are conjugate.     

    \begin{table}
      \begin{center} 
       \begin{tabular}{| c | c | c | c |}
    \hline
    Type & Class representatives & No. of Classes & No. of elements in a Class \\ \hline
    Central & $\begin{bmatrix} x & 0\\ 0 & x \end{bmatrix}$, $x\ne 0$ & $p-1$ & $1$ \\ \hline 
    Elliptic & $\begin{bmatrix} x & y\epsilon  \\ y & x \end{bmatrix}$, $y \ne 0$ &$ p(p-1)/2$ & $p^2-p$ \\
        \hline
        Parabolic &$\begin{bmatrix} x & 1\\ 0 & x \end{bmatrix}$, $x \ne 0$ & $p-1$ & $p^2-1$ \\ \hline
        Hyperbolic & $\begin{bmatrix} x & 0 \\ 0 & y \end{bmatrix}$,  $x \ne y$ & $(p-1)(p-2)/2$ & $p^2+p$ \\ \hline
    \end{tabular}
    \caption{Elements in $GL_2(F_p)$.  
    $x,y \in F_p$, $\epsilon \in F_p^\ast$ a non-square}\label{Ta:classes}
    \end{center}
    \end{table}
    
      The {\it central indicator} $d'$ of an element $\sigma  \in GL_2(F_p)$ is the smallest positive power of $\sigma$ which is contained in the center of $GL_2(F_p)$.   Equivalently, $d'$ is the order of $\Tilde{\sigma}$ in $PGL_2(F_p)$.   Thus $d'$ is a divisor of $d$, the order of $\sigma$, and we put $d=d'd''$, where $d''$ is the {\it central quotient}.  An {\it eigen-order} of an element is the multiplicative order of an eigenvalue (if there is one) in $F_p^\ast$.  

Let
\begin{equation}\label{E:classTypes} 
_pC(d', d''),\quad_pE(d', d''),\quad_pP^e (d', d''),\quad_pH^e_f (d', d'')
\end{equation}
denote the sets of central, elliptic, parabolic, hyperbolic elements, respectively, with central indicator $d'$,  central quotient $d'' = d/d'$, order $d'd''=d$, and  eigen-orders (if any) $e, f$.   We  omit the initial $p$ when it is understood or clear from the context.  We call the sets \eqref{E:classTypes} {\it class types} since they comprise smaller collections of conjugacy classes within the larger geometric types.  For hyperbolic elements, up to conjugacy, we may assume $e \leq f$.  For central elements, $d'=1$, and for parabolic elements, $d'=p$.    In the following lemmas, we give necessary and sufficient conditions (on the various parameters) for these sets to be non-empty, and, when non-empty, the number of conjugacy classes they contain.  Complete results for $p=3, 5, 7$ are given in Table~\ref{Ta:elements}.    
$\phi(n)$ denotes Euler's totient function, defined, for a positive integer $n$, as the number of positive integers less than or equal to $n$ and relatively prime to $n$.  In a cyclic group of order $n$, $\phi(d)$ is the number of elements of order $d\mid n$.
\subsection{Central and parabolic elements}
\begin{lemma}\label{L:centParaClasses}
Class types $_pC(1, d'')$ and $_pP^{d''} (p, d'')$ are non-empty if and only if $d'' \mid p-1$, in which case, each  contains $\phi(d'')$ conjugacy classes. 
\end{lemma}
\begin{proof}  A parabolic element is conjugate to $[\begin{smallmatrix} x & 1\\ 0 & x \end{smallmatrix}]$. Its $p$th power is the central element $[\begin{smallmatrix} x & 0\\ 0 & x \end{smallmatrix}]$.  Hence there is a conjugacy class in $\sigma \in\  _pP^{d''}(p, d'')$ if and only if $x \in F_p^\ast$ has multiplicative order $d''\mid p-1$.  There are $\phi(d'')$ such elements $x \in F_p^\ast$.  
\end{proof}
 \subsection{Elliptic elements}
     The $p^2 - 1$ elements 
   \begin{equation}\label{E:ellipticElements} \biggl\{\begin{bmatrix} x & y\epsilon  \\ y & x \end{bmatrix}, \quad \text{$x$, $y$ not both $=0$}\biggr\}
     \end{equation}
           form a cyclic subgroup of order $p^2-1$ isomorphic to $F_{p^2}^\ast$.    
 By analogy with complex numbers, one can think of these elements as  $x + y \sqrt\epsilon \in F_{p^2}$  ($\epsilon$ being a non-square in $F_p$.)
       Each element with $y\ne 0$ is  conjugate  to its $p$th power, 
     $$\begin{bmatrix} x & y\epsilon \\ y & x \end{bmatrix}^p = \begin{bmatrix} x & -y\epsilon \\ -y & x \end{bmatrix}.
     $$
     Moreover, 
     $$\begin{bmatrix} x & y\epsilon \\ y & x \end{bmatrix}^{p +1} = \begin{bmatrix} x^2 - y^2\epsilon & 0 \\ 0 & x^2 - y^2\epsilon \end{bmatrix}. 
     $$
  It follows that  the central indicator of an elliptic element is a divisor $d' >1$ of $p+1$.       Thus the general elliptic element has order $d=d'd''$, with $d' \mid p+1$, $d'>1$, and $d''\mid p-1$.     
  
\begin{lemma}\label{L:ellipClasses}
For $d' > 1$,   $d'\mid p+1$, and $d'' \mid p-1$, the class type $_pE(d', d'')$ is non-empty if and only if 
           \begin{equation}\label{E:ellipticCondition}
  \text{\rm gcd}((p-1)/d'', d') = 1,
  \end{equation}
  in which case  $_pE(d', d'')$ contains $\phi (d'd'')/2$ conjugacy classes.
  \end{lemma}
  \begin{proof}
 $\text{gcd}((p-1)/d'', d') \leq \text{gcd}(p-1, p+1)=2.$   Suppose $\text{gcd}((p-1)/d'', d')=2$. 
 Then $d=d'd''$ is even.  If $\sigma$ is an elliptic element of order $d=2k$,  then $\sigma^k
 $ is the central element 
  $$\begin{bmatrix} \epsilon^{(p-1)/2} & 0 \\
  0 & \epsilon^{(p-1)/2}\end{bmatrix} = \begin{bmatrix}-1 & 0 \\ 0& -1 \end{bmatrix}.$$
  This implies that $d''$ is even.  With both $d'$ and $d''$ even, 
  $d=4l$, and $\sigma^l$ is conjugate to one of the order $4$ elements
  $$\begin{bmatrix} 0 & \pm \epsilon^{(p+1)/4} \\
  \pm \epsilon^{(p-3)/4} & 0 \end{bmatrix},$$
  which implies $p \equiv -1 \mod 4$.
   On the other hand, the assumption $\text{gcd}((p-1)/d'', d')=2$ implies $(p-1)/d''$ is even.  But if $d''$ and $(p-1)/d''$ are both even,  $p \equiv 1 \pmod 4$, a contradiction.  Thus $\text{gcd}((p-1)/d'', d') = 1$ is a necessary condition.  
  
  To prove sufficiency,  note that every elliptic element is a power of an elliptic element $\sigma$ of maximum order $p^2-1$, with central indicator $p+1$.    Let $d'$ and  $d''$ satisfy the conditions of the lemma.   Let 
  $$n = \frac{p+1}{d'}\cdot \frac{p-1}{d''}.$$
  Then $\sigma^n$ has order $d=d'd''$ and $\sigma ^{nd'}$ is a central element of order $d''$.  If $d'$ is not the central indicator, there exists $d_1' \mid d'$, $1<d_1' < d'$, such that $\sigma^{nd_1'}$ is a central element of order
  $$ \dfrac{p-1}{\text{gcd}(\frac{p-1}{d_1'd''}, p-1)} = d_1'd''.
  $$
  Since $d_1'\mid \text{gcd}(p+1, p-1)$, the only possibility is $d_1' = 2$.  Then $d'$ is even,  and  condition \eqref{E:ellipticCondition} implies that $\frac{p-1}{d''}$ is odd, so that $\frac{p-1}{d_1'd''}$ is not an integer.  Thus  $d_1' = d'$ and $d'$ is the central indicator.
      \end{proof}

\subsection{Hyperbolic elements}

Let $x,y,z \geq 1$ be three positive integers with $\text{lcm}(x,y,z)=d$.  We say that $\{x,y,z\}$ satisfies the {\it strong lcm condition} if
$$\text{lcm}(x,y,z)=\text{lcm}(x,y) = \text{lcm}(x,z) = \text{lcm}(y,z).$$
An equivalent condition is that $\{x,y,z\}$ separates the prime divisors $q\mid d$ into two disjoint sets: 
\begin{enumerate}
\item[$U=$] primes $q\mid d$ having uniform multiplicity $m_q$ in the prime factorizations of $x,y,z$;
\item[$\overline U=$] primes $q\mid d$ whose multiplicity $j_q$ is non-maximal ($j_q < m_q$) in the prime factorization of exactly one of $x,y,z$.
\end{enumerate}

\begin{lemma}\label{L:hypClasses}
For $e,f,d' \mid p-1$, $e \leq f$, and $\text{lcm}(e,f,d')=d$, the class type  
$\ _pH^{e}_{f}(d', d/d')$ is non-empty if and only if  
$\{e,f,d'\}$ satisfies the strong lcm condition and the additional condition $2 \in \overline U$ if $2 \mid d$.  In this case, the number of conjugacy classes in the class type is 
\begin{equation}\label{E:hypClasses}
\prod_{q \in U} \phi(q^{m_q})q^{m_q-1}(q-2) \prod_{q \in \overline U} \phi(q^{m_q})\phi(q^{j_q})
\end{equation} 
if $e < f$, or half this number, if $e=f$.
\end{lemma}
\begin{proof}  First suppose $d = q^m$, a power of single prime.  A hyperbolic element of order $q^m$ is conjugate to a diagonal matrix $[\begin{smallmatrix} a & 0 \\ 0 & b \end{smallmatrix}]$, where $a,b  \in F_p^\ast$ have multiplicative orders $e=q^n$, $f=q^m$, and we assume $n \leq m$.  Let the central indicator $d'=q^j$, for some $j\leq m$.  Let $\text{D}_q (n,m,j)$ be the set of diagonal matrices with $(e,f,d')=(q^n, q^m, q^j)$.   We prove that $n<m$ implies $j=m$, and conversely, $j<m$ implies $n=m$ (i.e., that $q \in U$ or $q \in \overline U$).    

Associate diagonal matrices in $D_q(n,m,j)$ with pairs $(u,v) \in N \times N$, 
     $$N= \{0, 1, 2, \dots, q^m-1\},$$ 
     by choosing a generator $\gamma$ of the multiplicative group $F_{q^m}^\ast$, and defining the correspondence
$$(u,v) \leftrightarrow \begin{bmatrix} \gamma^u & 0 \\ 0& \gamma^v \end{bmatrix}.$$
 Membership in $\text{D}_q(n,m,j)$ is equivalent to the conditions 
\begin{equation}\label{E:pairCond}
\text{gcd}(u, q^m) = q^{m-n}, \quad \text{gcd}(v, q^m) = 1, 
\end{equation}
and, since $d'$ is the smallest solution of the congruence $ud' \equiv vd'$ mod $q^m$,
\begin{equation}\label{E:centIndCond}
\text{gcd}(u-v, q^m) = q^{m-j}.
\end{equation}
The set $\{u \in N \mid \text{gcd}(u, q^m) = q^{m-n}\}$ can be parametrized as
$$\{q^{m-n}(Aq + B) \mid A = 0, 1, \dots, q^{n-1}-1, \ B = 1, \dots, q-1\},$$
and the set $\{v \in N \mid \text{gcd}(v, q^m) = 1 \}$ as
$$\{aq+b \mid a = 0, 1, \dots, q^{m-1}-1, \ b = 1, \dots, q-1\}.$$ 
The differences $u-v$ between elements of the first and second sets can be laid out in an array of dimensions $\phi(q^n) \times \phi(q^m)$, 
$$\Delta = \{q^{m-n}(Aq+B) - (aq +b)\} = \{q^{m-n}[q(A-a) + B] - b\}.$$
If $n < m$, $\Delta$ consists of elements  $\equiv -b \pmod q$, that is, $\not\equiv 0 \pmod q$.  Hence for $n<m$ it follows from condition \eqref{E:centIndCond} that $j=m$ is the only possibility.  
\begin{equation}\label{E:diagSubCount}\notag
 |\text{D}_q(n,m,j)| = 
\begin{cases} 0 &\text{if $0 \leq j \leq m-1$} \\
\phi(q^n) \phi(q^m)& \text{if $j=m$.} 
\end{cases} \quad (n < m)
\end{equation}
 
If $n=m$, the $\phi(q^m) \times \phi(q^m)$ array of differences is
$$\Delta = \{q(A-a) + B - b\}.$$
In this array, $u-v \not\equiv 0 \pmod q$ if and only if $B-b \not\equiv 0 \pmod q$, and this occurs in 
$$\dfrac{q-2}{q-1}(\phi(q^m))^2 = \phi(q^m)q^{m-1}(q-2)
$$    
positions.  The remaining $\dfrac{1}{q-1}(\phi(m))^2$ positions have $u-v = q(A-a)\equiv q^{m-j} \pmod {q^m},$ $1 \leq j \leq m$, or, equivalently,
$$A - a \equiv q^{m-1-j} \pmod{q^{m-1}}, \quad 0 \leq j \leq m-1.$$
For each $j = 0, 1, \dots, m-1$ there are $\phi(q^j)$ differences $\equiv q^{m-1-j}$ in each of the $\phi (q^m)$ rows of the array.    
Thus if $n=m$,  
\begin{equation}\notag
 |\text{D}_q(m,m,j)| = 
\begin{cases} \phi(q^j)\phi(q^m) &\text{if $0 \leq j \leq m-1$} \\
\phi(q^m)q^{m-1}(q-2)& \text{if $j=m$.} 
\end{cases}
\end{equation}

We assumed $n\leq m$ in the triples $(n,m,j)$, but when $n \ne m$, 
$$|\text{D}_q(n,m,m)| = |\text{D}_q(m,n,m)|,$$
since conjugacy gives a one-to-one correspondence between the two sets.
Hence  we have shown that 
\begin{equation}\label{E:diagCount}
|\text{D}_q(j,m,m)|
=|\text{D}_q(m,j,m)|
 = |\text{D}_q(m,m,j)| =
 \begin{cases} \phi(q^j)\phi(q^m) & \text{if $j < m$} \\
 \phi(q^m)q^{m-1}(q-2) &\text{if $j=m$,}
 \end{cases}
 \end{equation}
 while, for all other triples $(i,j,k)$, $ 0 \leq i,j,k \leq m$, $|\text{D}_q(i,j,k)| = 0$.  
   
   The general hyperbolic element of order $d=\text{lcm}(x,y,z)$ is conjugate to a  product of diagonal matrices of prime power order.   The formula \eqref{E:hypClasses} is the product of \eqref{E:diagCount} over the set $U \cup \overline U$ of prime divisors of $d$.  The necessity of the $U, \overline U$ bipartition of the prime divisors of $d$, equivalently, of the strong lcm condition on $\{e,f,d'\}$, follows immediately.     The formula would yield $0$ if $2 \in U$, due to the factor $q-2$, which explains the extra condition.  
  
        If $e=f$  (and $d'   \ne 1$), the diagonal matrices $[\begin{smallmatrix} a & 0 \\ 0 & b \end{smallmatrix}]$ and $[\begin{smallmatrix} b & 0 \\ 0 & a \end{smallmatrix}]$ are distinct but conjugate.  This explains the necessity of reducing the formula by a factor of $\frac 12$ in this case.
\end{proof}

\begin{table}[h]
      \begin{center} 
       \begin{tabular}{| c || c | c || c | c || c | c || c | c ||  }
       \hline
     & \multicolumn{2}{| c ||}{Central} & \multicolumn{2}{| c ||}{Elliptic} &\multicolumn{2}{| c ||}{Parabolic} &\multicolumn{2}{| c ||}{Hyperbolic} 
      \\ \hline
     {\bf p}&$C(1, d'')$ & No. classes & $E(d', d'')$ & No. classes & $P^e(p, d'')$ & No. classes & $H^e_f(d', d'')$ & No. classes
      \\ \hline\hline
        {\bf 3} &  $C(1,1)$ & $1$ & $E(2,2)$ & $1$ & $P^1(3, 1)$ & $1$  & $H^1_2(2,1)$ & $1$\\
        \hline
        & $C(1,2)$ & $1$ & $E(4,2)$ & $2$  & $P^2(3, 2)$ & $1$ & & \\
        \hline 
        \\
        \hline
        {\bf 5} &$C(1,1)$ & $1$ & $E(2,4)$  & $2$ & $P^1(5,1)$ & $1$&$H^4_4(2,2)$  & $1$ \\ \hline
         &$C(1,2)$ & $1$ &$E(3,1)$ & $1$ & $P^2(5,2)$ & $1$ &$H^2_4(4,1)$ & $2$\\ \hline
          &$C(1,4)$ & $2$ &$E(3,2)$ &$1$  &$P^4(5,4)$ & $2$ & $H^1_4(4,1)$ &$2$ \\ \hline
           & &  & $E(3,4)$ & $2$ & &  &$H^1_2(2,1)$ &$1$ \\ \hline
           & &  &$E(6,4)$ & $4$ & &  & & \\ 
           \hline \\
           \hline
           {\bf 7} &$C(1,1)$ & $1$ &$E(2,2)$ & $1$ &$P^1(7,1)$ &$1$ &$H^6_6(3,2)$ &$1$ \\ \hline
           &$C(1,2)$ & $1$ & $E(2,6)$ & $2$ &$P^2(7, 2)$ & $1$ &$H^3_6(6,1)$ & $2$\\ \hline
           &$C(1,3)$ &$2$  & $E(4,2)$ & $2$ &$P^3(7,3)$ & $2$ &$H^3_6(2,3)$ & $2$ \\ \hline
           & $C(1,6)$ & $2$ &$E(4,6)$ & $4$ &$P^6(7,6)$ & $2$ &$H^2_6(3,2)$ &$2$ \\ \hline
           & &  & $E(8,2)$ &$4$  & &  & $H^1_6(6,1)$ &$2$\\ \hline
           & &  & $E(8,6)$ & $8$  & &  &$H^2_3(6,1)$ & $2$\\ \hline
           & & & & & & & $H^3_3(3,1)$ & $1$ \\ \hline
           & &  & &  & &  &$H^1_3(3,1)$ &$2$ \\ \hline
           & &  & &  & &  &$H^1_2(2,1)$ & $1$ \\ \hline
        \end{tabular}
        \caption{Class types  in $GL_2(F_p)$, $p=3,5,7$.}\label{Ta:elements}
    \end{center}
         \end{table}
   \section{The cardinality of $\text{Fix}_{\mathcal P}(\sigma)$}\label{S:action}
    The cardinality of $\text{Fix}_{\mathcal P}(\sigma) = \{ M \in{\mathcal M}({\mathcal P}) \mid \sigma M = M , \sigma \in GL_2(F_p)\}$ depends only on the conjugacy class of $\sigma$.  This follows from:  (i) $M \in {\mathcal M}({\mathcal P})$ if and only if $\gamma M \in {\mathcal M}({\mathcal P})$, $\gamma \in GL_2(F_p)$; and (ii) $\sigma M = M$ if and only if $\gamma\sigma\gamma^{-1}(\gamma M) = \gamma M$.          
    If $M \in \text{Fix}_{\mathcal P}(\sigma)$ then $M$ is a concatenation of $\sigma$-orbits on $V_2(F_p)$.  We do not say `disjoint union' because $M$, being a multiset, may contain the same orbit more than once. 
 \subsection{The zero-sum condition}
 \begin{lemma}\label{L:zSum} 
 If $\sigma \in GL_2(F_p)$ has no eigenvalue equal to $1$, 
 $\sigma M = M$ implies the zero-sum condition on $M$.
 \end{lemma}
 \begin{proof}  Let $\sigma$ be an arbitrary element of order $d > 1$.  $\sigma^d - I_2 = [\begin{smallmatrix} 0 & 0 \\ 0 & 0 \end{smallmatrix}]$, and if there is no eigenvalue equal to $1$, $\sigma - I_2$ is invertible. Therefore,
 \begin{equation}\label{E:matrixCongruence}
 (\sigma^d - I_2)(\sigma - I_2)^{-1} \equiv I_2 + \sigma + \sigma^2 + \dots + \sigma^{d-1} \equiv [\begin{smallmatrix} 0 & 0 \\ 0 & 0 \end{smallmatrix}]  \pmod p.
 \end{equation}
 Hence the $\sigma$-orbit of any column $[\begin{smallmatrix} a \\ b \end{smallmatrix}]$  of $M$  
  sums to 
  $$\begin{bmatrix} a \\b \end{bmatrix} (I_2 + \sigma +\sigma^2 + \dots + \sigma^{d-1}) \equiv \begin{bmatrix} 0 \\0 \end{bmatrix} \pmod p.$$
\end{proof}
   If $\sigma$ has an eigenvalue equal to $1$, it belongs to one of the conjugacy classes in $\ _pP^1(p, 1)$ or $ _pH^1_f(f, 1)$.     By suitable change of basis, we may take $\sigma$ to be 
  \begin{equation}\label{E:OneEigenElts}
  \begin{bmatrix} 1 & 1 \\ 0 & 1 \end{bmatrix} \in\  _pP^1(p, 1) \quad\text{or}\quad \begin{bmatrix} 1 & 0 \\
0 & y \end{bmatrix} \in\   _pH^1_f(f, 1).
\end{equation}
$\overline\infty$ is the $1$-eigenspace in both of these cases, and $\overline 0$ is the $f$-eigenspace in the hyperbolic case.  
\begin{lemma}\label{L:zSum2} Let  $M$ be a multiset in the form \eqref{E:pgenmatrix}.  If $\sigma M = M$,  $\sigma \in \  _pP^1(p, 1)$, the zero-sum condition is 
  \begin{equation}\label{E:pbolicZsum}
  \sum S_\infty \equiv \begin{bmatrix} 0 \\ 0 \end{bmatrix} \pmod p.
  \end{equation}
If $\sigma M = M$ and $\sigma \in \ _pH^1_f(f, 1)$, the zero-sum condition is 
\begin{equation}\label{E:hbolicZsum}
    \sum S_\infty \equiv -\sum \begin{bmatrix} S_{1}    \mid &   \dots \mid &  S_{p-1} \end{bmatrix} \pmod p,
    \end{equation}
    and in addition, 
    \begin{equation}\label{E:freeParam}
     \sum S_0 \equiv \begin{bmatrix} 0 \\0
    \end{bmatrix} \pmod p.
    \end{equation}
\end{lemma}
\begin{proof}
   For $\sigma = [\begin{smallmatrix} 1 & 1 \\ 0 & 1 \end{smallmatrix}] \in\  _pP^1(p, 1)$, a column $[\begin{smallmatrix} a \\ b \end{smallmatrix}] \in M$ with $b \not\equiv 0$ (a non-eigenvector) belongs to the $\sigma$-orbit
  $$\begin{bmatrix} a \\ b \end{bmatrix}, \begin{bmatrix} a+b \\ b \end{bmatrix}, \begin{bmatrix} a+2b \\ b \end{bmatrix}, \dots, \begin{bmatrix} a+(p-1)b \\ b \end{bmatrix}$$
 of length $p$, which spans all subspaces of $PG(1, F_p)$ except $\overline\infty$ and sums to 
  \begin{equation}\label{E:pbolicGenSum}
  \begin{bmatrix} ap + bp(p-1)/2 \\ bp \end{bmatrix} \equiv \begin{bmatrix} 0 \\ 0 \end{bmatrix} \pmod p.
  \end{equation}
  Hence the only non-trivial condition is \eqref{E:pbolicZsum}.
    
For $\sigma  = [\begin{smallmatrix} 1 & 0 \\
0 & y \end{smallmatrix}]\in \ _pH^1_f(f, 1)$, a column $[\begin{smallmatrix} a \\ b \end{smallmatrix}] \in M$   belongs to the $\sigma$-orbit
    $$\begin{bmatrix} a \\ b \end{bmatrix}, \begin{bmatrix} a \\ yb \end{bmatrix}, \begin{bmatrix} a \\ y^2b \end{bmatrix}, \dots, \begin{bmatrix} a \\ y^{f - 1} b \end{bmatrix}$$
    which sums to
    \begin{equation}\label{E:hOrbitSum}
    \begin{bmatrix} f a \\
    b (1 + y +y^2 + \dots + y^{f -1})
    \end{bmatrix} \equiv \begin{bmatrix} f a \\ 0 \end{bmatrix} \pmod p,
    \end{equation}
        an element of $\overline\infty$.   This implies \eqref{E:hbolicZsum}.   Taking $a \equiv 0$, we obtain  \eqref{E:freeParam}.
    \end{proof}
      Note that both non-trivial zero sum conditions \eqref{E:pbolicZsum} and \eqref{E:hbolicZsum} are equivalent to  weighted multi-partition congruences involving only the entries in the top row of $M$.

 \subsection{The shape of a $\sigma$-orbit}

 Let $\sigma$ be an arbitrary element of $GL_2(F_p)$, and let $\Tilde{\sigma}$ denote the image of $\sigma$ in the projective group $PGL_2(F_p)$.  $\Tilde{\sigma}$  permutes the $p+1$ points of $PG(1, F_p)$.  A {\it proper} orbit of $\Tilde{\sigma}$ consists of $d'$ points, where $d'$ is the central indicator.   The number of proper orbits of $\Tilde{\sigma}$ is  
\begin{equation}\label{E:Qdef}
o(\sigma) = \begin{cases}
 \frac{p+1}{d'} & \text{if $\sigma$ is elliptic or central;}\\
 \frac{p-1}{d'} & \text{if $\sigma$ is hyperbolic;}\\
 1 & \text{if $\sigma$ is parabolic.}
 \end{cases}
 \end{equation}
 We define a {\it generic} orbit of $\sigma$ to be an orbit that projects to a proper orbit of $\Tilde{\sigma}$ under the projection $V_2(F_p) \rightarrow PG(1, F_p)$.  
      
       Let ${\bf v} \in V_2(F_p)$ belong to a generic orbit of $\sigma$, and let $\overline v$ be its image in $PG(1, F_p)$.    Then ${\bf v}, \sigma({\bf v}), \dots, \sigma^{d'-1}({\bf v})$ project onto the proper $\Tilde{\sigma}$-orbit $\overline v, \Tilde{\sigma}(\overline v), \dots, \Tilde{\sigma}^{d'-1}(\overline v)$ on $PG(1, F_p)$.   $\sigma^{d'}$ is a central element $[\begin{smallmatrix} y & 0 \\ 0 & y \end{smallmatrix}]$, $y \in F_p^\ast$, of multiplicative order $d''$.  The $\sigma$-orbit of ${\bf v}$ can be displayed in the shape of a $d' \times d''$ array
   $$X({\bf v},y) = \begin{bmatrix}
     {\bf v} &  y{\bf v} & \dots & y^{d''-1}{\bf v} \\
    {\sigma}({\bf v}) &  y{\sigma}({\bf v}) & \dots & y^{d''-1}\sigma({\bf v}) \\
      \vdots \\
      {\sigma}^{d'-1}({\bf v}) &  y{\sigma}^{d'-1}({\bf v}) & \dots & y^{d''-1}\sigma^{d'-1}({\bf v})
    \end{bmatrix}.
    $$
   $X({\bf v},y)$ is {\it marked} in the vertical (projective) dimension by the $\Tilde{\sigma}$-orbit
 $$\langle \Tilde{\sigma} \rangle (\overline v) = \{\overline v, \Tilde{\sigma}(\overline v), \dots, \Tilde{\sigma}^{d'-1}(\overline v)\},
 $$
  and in the horizontal (central) dimension by the $y$-orbit
  $$
  \langle y \rangle = \{1, y, y^2, \dots, y^{d''-1}\} \subseteq F_p^\ast.
  $$
  $X({\bf v}, y)$ is generated by any one of its entries, e.g., starting with $y^i{\bf v}$ rather than ${\bf v}$ yields the same entries, with the columns shifted cyclically.  However, starting with $z{\bf v}$, $z \in F_p^\ast$, $z \not\in \langle y \rangle$,  produces a different array, $X(z{\bf v}, y)$, disjoint from $X({\bf v}, y)$, with the same vertical marking but a different horizontal marking, namely
  $$z\langle y \rangle = \{z, zy, zy^2, \dots, zy^{d''-1}\} \subseteq F_p^\ast.$$
    The number of distinct horizontal markings is the number of distinct orbits of $\langle y \rangle$ on $F_p^\ast$, namely, $\frac{p-1}{d''}$. 
    
    If ${\bf v}$ were an eigenvector of $\sigma$, the array would be an {\it eigen-orbit} of shape $1 \times e$ (where $e$ is the order of the eigenvalue $y$), marked in the vertical dimension  by the singleton orbit $\{\Tilde{\sigma}(\overline v)\} = \{\overline v \}$ and in the horizontal dimension by the $y$-orbit $\langle y \rangle \subseteq F_p^\ast$.  If ${\bf v}$ were replaced by $z{\bf v}$, $z \in F_p^\ast$, $z \not\in \langle y \rangle$, we would obtain a different linear array $X(z{\bf v}, y)$, disjoint from $X({\bf v}, y)$ with the same vertical marking $\{\overline v\}$ and a different horizontal marking.  The number of distinct horizontal markings is the number of orbits of $\langle y \rangle$ on $F_p^\ast$, namely $\frac{p-1}{e}$.

It follows that, for each $\sigma \in GL_2(F_p)$,  the set $V_2(F_p) - \{[\begin{smallmatrix} 0 \\ 0 \end{smallmatrix}]\}$ is partitioned into marked $\sigma$-orbits of generic shape $d'\times d''$, together, possibly, with marked eigen-orbits of shape $1 \times e$ and/or $1 \times f$.
  
\subsection{Tiling the Ferrers diagram}\label{S:tilings}
If a partition can be expressed in the form 
\begin{equation}\label{E:FeasPart}
{\mathcal P}=(le, nf, (\pi_1d'')^{[m_1d']}, (\pi_2d'')^{[m_2d']}, \dots, (\pi_ud'')^{[m_ud']}),
\end{equation}
where one or both of the first two parts ($le$, $nf$) may be $0$ or possibly  equal to one or two of the other parts,  there exists a {\it tiled} Ferrers diagram of ${\mathcal P}$, as pictured below.  If $l>0$, a row of the diagram, which we take to be the first, consists of $le$ dots and  can be covered by $l$ `tiles' of dimension $1\times e$.  Similarly, if $n>0$, another row, which we take to be the second, consists of $nf$ dots which can be covered by $n$ tiles of size $1 \times f$.  The remaining rows consist of rectangular arrays of dots, of size $m_id' \times \pi_id''$, which can be covered by tiles of size $d' \times d''$ in an $m_i \times \pi_i$  arrangement, $i = 1, \dots , u$.  
    \begin{equation}\label{E:genericTile}
\begin{array}{l l ll l}
 &\qquad  \begin{tikzpicture}
 \draw[ultra thick,black] (0,0) -- (.2,0); 
 \draw[ultra thick,black] (.25,0) -- (.45,0);
 \draw[ultra thick,black] (.50,0) -- (.70,0);
  \draw[ultra thick,black] (.75,0) -- (.95,0);
 \end{tikzpicture}    \ \cdots \  \begin{tikzpicture}
 \draw[ultra thick,black] (0,0) -- (.2,0);
 \end{tikzpicture}
 &\qquad l(1 \times e) \\
 &\qquad  \begin{tikzpicture}
 \draw[ultra thick,black] (0,0) -- (.35,0);
 \draw[ultra thick,black] (.40,0) -- (.75,0);
 \draw[ultra thick,black] (.80,0) -- (1.15,0);
 \end{tikzpicture}  \ \cdots \ \begin{tikzpicture}
 \draw[ultra thick,black] (0,0) -- (.35,0);
 \end{tikzpicture}
& \qquad n(1 \times f)\\
\\
  &\qquad  \begin{tikzpicture}
  \draw[thick] (0,0) rectangle (.15, .3);
  \draw[thick] (.20,0) rectangle (.35, .3);
  \draw[thick] (.40,0) rectangle (.55, .3);
  \draw[thick] (.60,0) rectangle (.75, .3);
 \end{tikzpicture}
\ \cdots \ \begin{tikzpicture}
  \draw[thick] (0,0) rectangle (.15, .3);
   \end{tikzpicture}
&\qquad \\
&\qquad  \begin{tikzpicture}
  \draw[thick] (0,0) rectangle (.15, .3);
  \draw[thick] (.20,0) rectangle (.35, .3);
  \draw[thick] (.40,0) rectangle (.55, .3);
  \draw[thick] (.60,0) rectangle (.75, .3);
 \end{tikzpicture}
\ \cdots \ \begin{tikzpicture}
  \draw[thick] (0,0) rectangle (.15, .3);
   \end{tikzpicture}
    & \qquad   m_1d' \times \pi_1d'' \\
&\qquad  \vdots   \\
&\qquad  \begin{tikzpicture}
  \draw[thick] (0,0) rectangle (.15, .3);
  \draw[thick] (.20,0) rectangle (.35, .3);
  \draw[thick] (.40,0) rectangle (.55, .3);
  \draw[thick] (.60,0) rectangle (.75, .3);
 \end{tikzpicture}
\ \cdots \ \begin{tikzpicture}
  \draw[thick] (0,0) rectangle (.15, .3);
   \end{tikzpicture}
&\qquad \\
&\qquad \\ 

& \qquad \begin{tikzpicture}
  \draw[thick] (0,0) rectangle (.15, .3);
  \draw[thick] (.20,0) rectangle (.35, .3);
  \draw[thick] (.40,0) rectangle (.55, .3);
 \end{tikzpicture}
\ \cdots \ \begin{tikzpicture}
  \draw[thick] (0,0) rectangle (.15, .3);
   \end{tikzpicture}
 & \qquad \\
& \qquad \begin{tikzpicture}
  \draw[thick] (0,0) rectangle (.15, .3);
  \draw[thick] (.20,0) rectangle (.35, .3);
  \draw[thick] (.40,0) rectangle (.55, .3);
 \end{tikzpicture}
\ \cdots \ \begin{tikzpicture}
  \draw[thick] (0,0) rectangle (.15, .3);
   \end{tikzpicture}
 & \qquad m_2d' \times \pi_2d'' \\
 &\qquad  \vdots  & \\
& \qquad \begin{tikzpicture}
  \draw[thick] (0,0) rectangle (.15, .3);
  \draw[thick] (.20,0) rectangle (.35, .3);
  \draw[thick] (.40,0) rectangle (.55, .3);
 \end{tikzpicture}
\ \cdots \ \begin{tikzpicture}
  \draw[thick] (0,0) rectangle (.15, .3);
   \end{tikzpicture}
 & \qquad \\
    & & \qquad\vdots \\
   &  \\
 
& \qquad
 \begin{tikzpicture}
  \draw[thick] (0,0) rectangle (.15, .3);
 \end{tikzpicture}
 \ \cdots\ \begin{tikzpicture}
  \draw[thick] (0,0) rectangle (.15, .3);
   \end{tikzpicture}
   & \qquad \\
 &\qquad  \vdots  &\qquad m_ud' \times \pi_ud''  \\
& \qquad
 \begin{tikzpicture}
  \draw[thick] (0,0) rectangle (.15, .3);
 \end{tikzpicture}
\ \cdots\ \begin{tikzpicture}
  \draw[thick] (0,0) rectangle (.15, .3);
   \end{tikzpicture}
  \end{array}
 \qquad\qquad
 \begin{tabular}{l | l}
 \text{Tile} & \text{Dimensions} \\
 \hline 
 \begin{tikzpicture}
 \draw[ultra thick,black] (0,0) -- (.2,0);
 \end{tikzpicture} & $1 \times e$ \\
 \begin{tikzpicture}
 \draw[ultra thick,black] (0,0) -- (.35,0);
 \end{tikzpicture} & $1 \times f$  \\
 \begin{tikzpicture}
  \draw[thick] (0,0) rectangle (.15, .3);
   \end{tikzpicture} & $d' \times d''$
 \end{tabular}
\end{equation}

Such a diagram is  {\it feasible} for tiling by orbits of $\sigma \in \ _pH^e_f(d',d'')$.  If each $d' \times d''$ tile  is substituted by a marked generic $\sigma$-orbit, and each $1\times e$, $1 \times f$ tile  by a marked eigen-orbit, the result is a column multiset $M$ for which $\sigma M = M$.  If $\sigma \in GL_2(F_p)$  belongs to a different class type (central, parabolic, elliptic), a partition of the form \eqref{E:FeasPart} may still be feasible for tiling by $\sigma$-orbits, if certain additional restrictions are satisfied. 
  
  \begin{lemma}\label{L:NEW_feas}  A partition ${\mathcal P}$, expressible in the form \eqref{E:FeasPart}, is feasible for tiling by $\sigma$-orbits, $\sigma \in GL_2(F_p)$, if 
$\sigma$ has central indicator $d'$, central quotient $d''$, $u \leq o(\sigma)$, and the following additional conditions, depending on the class type of $\sigma$,   are satisfied:
$$\begin{tabular}{c | c}
{\rm Class Type of\ } $\sigma$ & \text{\rm Conditions on ${\mathcal P}$} \\
\hline
\text{\rm Central} & $l=n=0$, $u >1$ \\
\text{\rm Elliptic} &$ l=n=0$, $u > 0$ \\
\text{\rm Parabolic} & $n=0$, $u=m_1=1$, ($l\ne 1$)*
 \\
\text{\rm Hyperbolic} & \text{\rm If}\  $u=0$, \text{\rm then}\ $l,n > 0$, ($l\ne 1$)*.
\end{tabular}
$$
* If $\sigma$ has an eigenvalue equal to $1$.
\end{lemma}
\begin{proof}
Central and elliptic elements have no eigenspaces, so $l=n=0$.    For central elements,  generic orbits have shape $1 \times d''$, so $u >1$ is necessary to satisfy the rank condition on the multiset $M$.  For parabolic elements, $u \geq 1$   is necessary to satisfy the rank condition. Hence $u=1$, and $m_1=1$ since a single generic orbit spans all subspaces of $PG(1, F_p)$ except the eigenspace.  For hyperbolic elements, if $u=0$, then both $l$ and $n$ must be non-zero to satisfy the rank condition.  The extra condition $l \ne 1$, pertaining to the eigenvalue $1$ cases, is explained in the proof of Theorem~\ref{T:NEW_OneEigen}.  
\end{proof}

 Feasibility of ${\mathcal P}$ for tiling by $\sigma$-orbits is necessary, but not always sufficient, for the existence of a zero-sum multiset $M \in \text{Fix}_{\mathcal P}(\sigma)$.   
  
\begin{theorem}\label{T:NEW_noOneEigen}
 If $\sigma \in GL_2(F_p)$ has no eigenvalue equal to $1$, and ${\mathcal P}$ is a partition expressible in the form \eqref{E:FeasPart}, feasible for tiling by $\sigma$-orbits, then  
  \begin{equation}\label{E:NEW_NoEigenOne}\notag
|\text{Fix}_{\mathcal P}(\sigma)| = O(\sigma, {\mathcal P}) \binom{l+\frac{p-1}{e}-1}{l}^* \binom{n + \frac{p-1}{f} - 1}{n}^*\prod_{i =1}^u\binom{\pi_i+\frac{p-1}{d''}-1}{\pi_i}^{m_i},
  \end{equation}
  where 
    \begin{equation}\notag
O(\sigma, {\mathcal P})= \dfrac{o(\sigma)!}{(o(\sigma)-u)! m_1! m_2! \dots m_u !},
\end{equation}
and $1 < e < f$ are the  orders of the eigenvalues, if any.  If an eigenvalue does not  exist, the  corresponding factor, marked with an asterisk, is set equal to $1$.    \end{theorem}    
 \begin{proof} 
 A  multiset $M$ such that $\sigma M = M$ is realized by any assignment of vertical and horizontal markings to \eqref{E:genericTile}.  The zero-sum condition on $M$ is guaranteed by the absence of an eigenvalue equal to $1$.    A  vertical (projective) marking  is an unordered selection of $m_i$ distinct proper $\Tilde{\sigma}$-orbits $\langle \Tilde{\sigma} \rangle \{\overline{v_{ij}}\}$, $j=1, \dots, m_i$, for each block of generic tiles, $i = 1, \dots, u$.  
  The number of distinct possible selections is $O(\sigma, {\mathcal P})$.
   Each of the $m_i\pi_i$ generic tiles in the $i$th block of  generic tiles can have any one of $\frac{p-1}{d''}$  horizontal markings.   Similarly, each tile in a (possible)  eigenspace can have any one of  $\frac{p-1}{e}$ or $\frac{p-1}{f}$ horizontal markings.    
 Within a row, tiles are unordered.  We can view them as `balls' placed arbitrarily into `boxes' labelled with the available markings.  The general balls-in-boxes coefficient \eqref{E:ballsInBoxes} accounts for the remaining factors in the formula.   If eigenspaces exist but are unmarked ($l$ or $n$ $=0$), the corresponding binomial coefficients reduce to $1$.  
\end{proof}  

\begin{theorem}\label{T:NEW_OneEigen} 
If $\sigma \in GL_2(F_p)$ has an eigenvalue equal to $1$, and ${\mathcal P}$ is a partition expressible in the form \eqref{E:FeasPart}, feasible for tiling by $\sigma$-orbits, then 
\begin{equation}\label{E:NEW_hypTileCount}
|\text{Fix}_{\mathcal P}(\sigma)| =
\begin{cases}
 O(\sigma, {\mathcal P})\displaystyle{\binom{n+\frac{p-1}{f} -1}{n}W_{(l, \pi_1^{[m_1]}, \dots, \ \pi_u^{[m_u]})}} &\text{if $\sigma \in \ _pH^1_f(f, 1)$,} \\
   \displaystyle{\binom{\pi_1+p-2}{\pi_1}W_{(l)}} &\text{if $\sigma \in \ _pP^1(p, 1)$.}
   \end{cases}
  \end{equation}
where $O(\sigma, {\mathcal P})$ has the definition given in Theorem~\ref{T:NEW_noOneEigen}, and  $W_{(-)}$ has the definition given in Theorem~\ref{T:multiPart}.

  \begin{proof}
  As before, a vertical  marking  is an unordered selection of $m_i$ distinct proper $\Tilde{\sigma}$-orbits $\langle \Tilde{\sigma} \rangle \{\overline{v_{ij}}\}$, $j=1, \dots, m_i,$ for each block of generic tiles, $i = 1, \dots, u$.  In the cases at hand, the number of possible vertical markings is 
  \begin{equation}\notag
  O(\sigma, {\mathcal P}) =
  \begin{cases} 
  \dfrac{(\frac{p-1}{f})!}{(\frac{p-1}{f} -u)! m_1! m_2! \dots m_u !} & \text{if $\sigma \in \ _pH^1_f(f, 1)$; and} \\
  1 &\text{if $\sigma \in \ _pP^1(p, 1)$.}
   \end{cases}
\end{equation}
In the hyperbolic case we may assume without loss of generality that $\sigma = [\begin{smallmatrix} 1 & 0 \\ 0 & y \end{smallmatrix}]$, with eigenspaces are $\overline\infty$ and $ \overline 0$.  Then the set $\{\overline{v_{11}}, \overline{v_{12}}, \dots, \overline{v_{1m_1}},\dots, \overline{v_{um_u}}\}$ consists of distinct vectors of the form $ [\begin{smallmatrix} a_{ij} \\ 1\end{smallmatrix}]$,  $1 \leq a_{ij} \leq p-1$.  The horizontal markings of the $1 \times 1$ tiles in the $1$-eigenspace and of the $f \times 1$ generic  tiles are just choices of scalars from $F_p^\ast$, which must collectively satisfy the zero-sum condition \eqref{E:hbolicZsum}.      Let $r_1, \dots, r_l$ be the scalar markings of the $1\times 1$ tiles, and $z_{ij1}, \dots, z_{ij\pi_i }$, $i = 1, \dots, u$, $j=1, \dots m_i$, be the scalar markings of tiles  in the $i$th block of generic tiles.  Then  \eqref{E:hbolicZsum} is a weighted multi-partition congruence, with tuple $(l, \pi_1 m_1, $ $\dots, \pi_u m_u)$ and weights $(1, fa_{11}, \dots, fa_{um_u})$.   
 \begin{equation}\label{E:New_hypMultPart}
r_1 + \dots +r_l  + \sum_{i=1}^u  \sum_{j=1}^{m_i} fa_{ij}\sum_{k=1}^{\pi_i}z_{ijk} \equiv 0 \pmod p.
\end{equation}
  The number of solutions,  by Theorem~\ref{T:multiPart}, is $W_{(l, \pi_1^{[m_1]}, \dots, \ \pi_u^{[m_u]})}$.  By the proof of the theorem, this number is independent of the weights.     The horizontal marking of the $f$-eigenspace can be chosen arbitrarily in
$$
\binom{n + \frac{p-1}{f} - 1}{n}
$$
ways, in view of \eqref{E:freeParam}.  If $u=0$, $W_{(l)} = 0$ if $l=1$.  Thus we have the additional feasibility condition $l\ne 1$ in this case.

In the parabolic case, we may assume without loss of generality that $\sigma =\begin{bmatrix} 1 & 1 \\ 0 & 1 \end{bmatrix}$,  with eigenspace $\overline\infty$.  The single row of $m>0$ generic tiles of shape $p \times 1$ can be projectively marked in only one way by the single proper orbit $\langle \Tilde{\sigma}\rangle(\overline 0)$ $=\{\overline 0, \overline 1, \dots, \overline{p-1}\}$.  Generic $\sigma$-orbits are zero-sum (cf. \eqref{E:pbolicGenSum}), so the generic tiles may be horizontally marked by scalars from $F_p^\ast$ arbitrarily, in 
$$\binom{m + p-2}{m}$$
ways,  by the balls-in-boxes analogy.   $1 \times 1$ tiles in the top (eigenspace) row are also marked by scalar choices from $F_p^\ast$.   The zero-sum condition \eqref{E:pbolicZsum} is simply that these  $l\geq 0$  scalars sum to $0$ mod $p$.  By Theorem~\ref{T:multiPart}, the number of scalar choices satisfying this condition is $W_{(l)}$.  Since $W_{(1)} = 0$,  we have the additional feasibility condition $l\ne 1$ in the parabolic case as well.

  \end{proof}
\end{theorem}

 \subsection{Multiple tilings} 
  
 If $\sigma \in \ _pH^e_f(d',d'')$, it is possible for a partition ${\mathcal P}$  to admit several distinct expressions of the form \eqref{E:FeasPart}, each feasible for tiling by orbits of $\sigma$.  

Let ${\mathcal P}$  be expressed in the form
$${\mathcal P}=(x,y, (\pi_1d'')^{[m_1d']}, (\pi_2d'')^{[m_2d']}, \dots, (\pi_ud'')^{[m_ud']} ),\quad x,y \geq 0.$$
\begin{enumerate}
\item If $x, y >0$, $x \ne y$, $e \ne f$, and both $x$ and $y$ are divisible by both $e$ and $f$, the eigenspaces can be marked in two distinct ways:  (i) by  $x/e$ tiles of shape $1 \times e$  in the $e$-eigenspace of $\sigma$, and by $y/f$ tiles of shape $f \times 1$  in the $f$-eigenspace; or (ii) by $y/e$ tiles of shape $1 \times e$  in the $e$-eigenspace and $x/f$ tiles of shape $1 \times f$  in the $f$-eigenspace.  If $1<e < f$, the product of the asterisked factors in Theorem~\ref{T:NEW_noOneEigen} is replaced by
$$
 \biggl[\binom{\frac xe +\frac{p-1}{e}-1}{\frac xe} \binom{\frac yf + \frac{p-1}{f} - 1}{\frac yf} +  
 \binom{\frac ye +\frac{p-1}{e}-1}{\frac ye} \binom{\frac xf + \frac{p-1}{f} - 1}{\frac xf} \biggr].
$$
\\

\item If $x>0$, $y=0$, $e \ne f$ and $x$ is divisible by both $e$ and $f$, either of the two eigenspaces can be marked:  (i) the $e$-eigenspace by $x/e$ tiles of shape $1 \times e$; or (ii) the $f$-eigenspace by $x/f$ tiles of shape $1 \times f$. If $1 < e < f$, the product of the asterisked factors in Theorem~\ref{T:NEW_noOneEigen} is replaced by
$$
 \biggl[\binom{\frac xe +\frac{p-1}{e}-1}{\frac xe} +  
 \binom{\frac xf + \frac{p-1}{f} - 1}{\frac xf} \biggr].
$$
 \\

\item If $x=y=0$ and $\sigma$ has  central indicator $d'=2$, there can be up to $u$ distinct tilings.    Hyperbolic  elements with central indicator $2$ (which are characterized as traceless) fall into two types, according to Lemma~\ref{L:hypClasses}:
\begin{align}\notag
_pH^e_{2e}(2, e), & \qquad\text{$e \equiv 1 \pmod 2$; or} \\ \notag
_pH^e_{e}(2,\frac e2), & \qquad\text{$e\equiv 0 \pmod 4$}.
\end{align}
If $e \equiv 1 \pmod 2$, and there is a part $\pi_j e \in {\mathcal P}$ divisible by $2e$, its multiplicity can be reduced from $2m_j$ to $2(m_j - 1)$, putting $x=y = \pi_j e$, and using those two parts to mark the $e$- and $2e$-eigenspaces.  If $e \equiv 0 \pmod 4$, and there is a part $\pi_j e/2$ divisible by $e$, its multiplicity can be reduced in the same way,  putting $x=y=\pi_je/2$, and using those two parts to mark the two $e$-eigenspaces.  If $e \equiv 1 \pmod p$, $e \ne 1$, the formula for $|\text{Fix}_{\mathcal P}(\sigma)|$ given by Theorem~\ref{T:NEW_noOneEigen},
$$O(\sigma, {\mathcal P}) \prod_{i =1}^u\binom{\pi_i+\frac{p-1}{d''}-1}{\pi_i}^{m_i},
$$
is increased by a sum of terms of the form
\begin{equation}\label{E:multipleTerms}
O_j(\sigma, {\mathcal P})\binom{\frac{\pi_j}{e}+\frac{p-1}{e}-1}{\frac{\pi_j}{e}} \binom{\frac{\pi_j}{2e} + \frac{p-1}{2e} - 1}{\frac{\pi_j}{2e}}\prod_{i =1}^u\binom{\pi_i+\frac{p-1}{d''}-1}{\pi_i}^{m_i^*},
\end{equation}
where
\begin{equation}\label{E:Oj}
O_j(\sigma, {\mathcal P}) = \dfrac{o(\sigma)!}{(o(\sigma)-u+1)! m_1^\ast! m_2^\ast! \dots  m_u^\ast !},\end{equation}
and
\begin{equation}\label{E:mast}
 m_i^\ast =
\begin{cases} m_i & \text{if $i \ne j$;}\\
m_i-1 & \text{if $i=j$.}
\end{cases}
\end{equation}
There is one such term for each part $\pi_j e$ divisible by $2e$.  If $e \equiv 0 \pmod 4$ the product
$$\binom{\frac{\pi_j}{e}+\frac{p-1}{e}-1}{\frac{\pi_j}{e}} \binom{\frac{\pi_j}{2e} + \frac{p-1}{2e} - 1}{\frac{\pi_j}{2e}}$$
 in \eqref{E:multipleTerms} is replaced by 
$$\binom{\frac{\pi_j}{e}+\frac{p-1}{e}-1}{\frac{\pi_j}{e}}^2,$$
and there is one  term for each part $\pi_j e/2$ divisible by $e$.
\end{enumerate}

If $e=1$, in the three cases above, adjustments to the formula for  $|\text{Fix}_{\mathcal P}(\sigma)|$ in Theorem~\ref{T:NEW_OneEigen} are as follows:
\begin{enumerate}
\item $O(\sigma, {\mathcal P})\biggl[\displaystyle{\binom{\frac xf+\frac{p-1}{f} -1}{\frac xf}} W_{(y, \pi_1^{[m_1]}, \dots, \ \pi_u^{[m_u]})} + \displaystyle{\binom{\frac yf+\frac{p-1}{f} -1}{\frac yf}} W_{(x, \pi_1^{[m_1]}, \dots, \ \pi_u^{[m_u]})}\biggr]$;
\\
\item $O(\sigma, {\mathcal P})\biggl[\displaystyle{\binom{\frac xf+\frac{p-1}{f} -1}{\frac xf}} W_{(\pi_1^{[m_1]}, \dots, \ \pi_u^{[m_u]})} +  W_{(x, \pi_1^{[m_1]}, \dots, \ \pi_u^{[m_u]})}\biggr]$;
\\
\item For $\sigma \in \ _pH^1_2(2, 1)$, the formula given in Theorem~\ref{T:NEW_OneEigen}, 
$$O(\sigma, {\mathcal P}) W_{(\pi_1^{[m_1]}, \dots, \ \pi_u^{[m_u]})},$$
is increased by a sum of terms of the form 
$$O_j(\sigma, {\mathcal P})\displaystyle{\binom{\frac{\pi_j}{2} + \frac{p-1}{2} - 1}{\frac{\pi_j}{2}}} W_{(\pi_j, \pi_1^{[m_1^*]}, \dots, \ \pi_u^{[m_u^*]})},$$
with $O_j(\sigma, {\mathcal P})$ and $m_i^\ast$ as defined above at \eqref{E:Oj} and \eqref{E:mast}.  There is one term for each part $\pi_j$ divisible by $2$.
\end{enumerate}

As an example of the last type, let ${\mathcal P} = (4^{[2]}, 2^{[4]}, 1^{[2]})$ and $\sigma \in\ _pH^1_2(2,1)$.  Then we have
\begin{equation}\notag
\begin{split}
|\text{Fix}_{\mathcal P}(\sigma)| &=O(\sigma, {\mathcal P}) W_{(4^{[2]}, 2^{[4]}, 1^{[2]})} \\
&\quad +O_1(\sigma, {\mathcal P})\frac{p^2-1}{8} W_{(4, 2^{[4]}, 1^{[2]} )} \\
&\qquad + O_2(\sigma, {\mathcal P})\frac{p-1}{2} W_{(2, 4^{[2]},2^{[2]}, 1^{[2]} )}.
\end{split}
\end{equation}
The three terms correspond to the three tilings 
    \begin{equation}\label{E:NEW_ex1}\notag
\begin{array}{l l l l l l l l l l l l l}
    {\mathcal T}_1 &= & \begin{tikzpicture}
 \draw[ultra thick,black] (0,1.4) -- (0,1.2);
  \draw[ultra thick,black] (.2,1.4) -- (.2,1.2);
  \draw[ultra thick,black] (.4,1.4) -- (.4,1.2);
  \draw[ultra thick,black] (.6,1.4) -- (.6,1.2);
\draw[ultra thick,black] (0,1) -- (0,.8);
  \draw[ultra thick,black] (.2,1) -- (.2,.8);
 \draw[ultra thick,black] (0,.6) -- (0,.4);
  \draw[ultra thick,black] (.2,.6) -- (.2,.4);
   \draw[ultra thick,black] (0,.2) -- (0,0);
 \end{tikzpicture}
  &  {\mathcal T}_2  &=&  \begin{tikzpicture}
  \filldraw[black] (0,1.4) circle (1.5pt);  
 \filldraw[black] (.2,1.4) circle (1.5pt);
  \filldraw[black] (.4,1.4) circle (1.5pt);  
 \filldraw[black] (.6,1.4) circle (1.5pt);
  \draw[ultra thick,black] (-.05,1.2) -- (.3,1.2);
  \draw[ultra thick,black] (.35,1.2) -- (.65,1.2);
\draw[ultra thick,black] (0,1) -- (0,.8);
  \draw[ultra thick,black] (.2,1) -- (.2,.8);
 \draw[ultra thick,black] (0,.6) -- (0,.4);
  \draw[ultra thick,black] (.2,.6) -- (.2,.4);
   \draw[ultra thick,black] (0,.2) -- (0,0);
   \end{tikzpicture}  
   &  {\mathcal T}_3  &=&  \begin{tikzpicture}
  \filldraw[black] (0,1.4) circle (1.5pt);  
 \filldraw[black] (.2,1.4) circle (1.5pt);
  \draw[ultra thick,black] (-.05,1.2) -- (.3,1.2);
\draw[ultra thick,black] (0,1) -- (0,.8);
  \draw[ultra thick,black] (.2,1) -- (.2,.8);
  \draw[ultra thick,black] (.4,1) -- (.4,.8);
  \draw[ultra thick,black] (.6,1) -- (.6,.8);
 \draw[ultra thick,black] (0,.6) -- (0,.4);
  \draw[ultra thick,black] (.2,.6) -- (.2,.4);
   \draw[ultra thick,black] (0,.2) -- (0,0);
   \end{tikzpicture}  
   &\quad \text{of}\quad
   & {\mathcal P} &=& \begin{tikzpicture}
   \filldraw[black] (0,1.4) circle (1.5pt);  
 \filldraw[black] (.2,1.4) circle (1.5pt);
  \filldraw[black] (.4,1.4) circle (1.5pt);  
 \filldraw[black] (.6,1.4) circle (1.5pt);
 \filldraw[black] (0,1.2) circle (1.5pt);  
 \filldraw[black] (.2,1.2) circle (1.5pt);
 \filldraw[black] (.4,1.2) circle (1.5pt);  
 \filldraw[black] (.6,1.2) circle (1.5pt);
      \filldraw[black] (0,1) circle (1.5pt);  
 \filldraw[black] (.2,1) circle (1.5pt);
    \filldraw[black] (0,.8) circle (1.5pt);  
 \filldraw[black] (.2,.8) circle (1.5pt);
    \filldraw[black] (0,.6) circle (1.5pt);  
 \filldraw[black] (.2,.6) circle (1.5pt);
  \filldraw[black] (0,.4) circle (1.5pt);  
 \filldraw[black] (.2,.4) circle (1.5pt);
  \filldraw[black] (0,0) circle (1.5pt);  
 \filldraw[black] (0,.2) circle (1.5pt);
  \end{tikzpicture}\quad .
\end{array}
\end{equation}
Vertical line segments denote $2 \times 1$ generic tiles; horizontal segments denote $1 \times 2$ tiles in the $2$-eigenspace.

\section{Burnside's Lemma}\label{S:Burnside}

Given ${\mathcal P}$, we may now replace the summation over the elements of $GL_2(F_p)$ in Burnside's Lemma \eqref{E:Burnside} by a sum over a subset consisting of the identity element and  one representative of each class type for which ${\mathcal P}$ has a feasible tiling.   The identity contributes the summand $|{\mathcal M}({\mathcal P})|$, obtained from Theorem~\ref{T:generalP} or its generalization.
   Each remaining summand is a three-fold product of the form 
\begin{equation}\label{E:threeFold}
t(\sigma)c(\sigma) |\text{Fix}_{\mathcal P}(\sigma)|,
\end{equation}
where $\sigma$ is a representative of a feasible class type; $t=t(\sigma)$ is the number of elements in the class (obtained from Table~\ref{Ta:classes}); and $c = c(\sigma)$ is the number of classes in the class type (obtained from Lemmas~\ref{L:centParaClasses}, ~\ref{L:ellipClasses},  or ~\ref{L:hypClasses}).  These factors are independent of the choice of $\sigma$ within the class type.  Let $F({\mathcal P})\subset GL_2(F_p)$ be a set containing one representative of each class type for which ${\mathcal P}$ has a feasible tiling.  Let $O({\mathcal P})$ be the set of orbits of $GL_2(F_p)$ on ${\mathcal M}({\mathcal P})$.  Burnside's Lemma now reads
\begin{equation}\label{E:specialBurnside} 
|O({\mathcal P})| = (p^2-1)^{-1}(p^2-p)^{-1}\bigl\{|{\mathcal M}({\mathcal P})| +  \sum_{\sigma \in F({\mathcal P})}t(\sigma)c(\sigma) \sum_{{\mathcal T}\supset {\mathcal P}}|\text{Fix}_{\mathcal T}(\sigma)|\bigr\}.\end{equation} 
We use the notation ${\mathcal T} \supset {\mathcal P}$ in the cases where ${\mathcal T}$ is one of a multiplicity of feasible tilings of ${\mathcal P}$ by $\sigma$-orbits.  In these cases $|\text{Fix}_{\mathcal P}(\sigma)|=\sum_{{\mathcal T}\supset {\mathcal P}}|\text{Fix}_{\mathcal T}(\sigma)|$.

The computations summarized in the tables below  yield $|O({\mathcal P})| = 2, 6$ for the $6$-partition ${\mathcal P}=(2^{[3]})$ and the primes $p=3,5$ respectively.
$$
\begin{array}{l l}
       \begin{tabular}{| c | c | c  |  }
       \hline
    ${\mathcal T}$ & Class Type & Summand  
      \\ \hline\hline
        \rule{0ex}{3.5ex}
        \begin{tikzpicture}
  \filldraw[black] (0,.4) circle (1.5pt);  
 \filldraw[black] (.2,.4) circle (1.5pt);
  \filldraw[black] (0,.2) circle (1.5pt);  
 \filldraw[black] (.2,.2) circle (1.5pt);
  \filldraw[black] (0,0) circle (1.5pt);  
 \filldraw[black] (.2,0) circle (1.5pt);
  \end{tikzpicture}
 & $Id$ & $12$     \\
 \hline
 \rule{0ex}{3.5ex}
  \begin{tikzpicture}
 \draw[ultra thick,black] (-.05,.5) -- (.25,.5);
 \draw[ultra thick,black] (-.05, .25) -- (.25,.25);
  \draw[ultra thick,black] (-.05,0) -- (.25,0);
   \end{tikzpicture}   &$\ _3C(1,2)$ &$4$  \\
   \hline
 \rule{0ex}{3.5ex}
 \begin{tikzpicture}
 \filldraw[black] (0,.5) circle (1.5pt);  
 \filldraw[black] (.2,.5) circle (1.5pt);
 \draw[ultra thick,black] (0,0) -- (0,.3);
  \draw[ultra thick,black] (.2,0) -- (.2,.3);
 \end{tikzpicture} & $\ _3H^1_2(2,1)$&$36$  \\
        \hline
        \rule{0ex}{3.5ex}
         \begin{tikzpicture}
 \draw[ultra thick,black] (-.05,.5) -- (.25,.5);
 \draw[ultra thick,black] (0,0) -- (0,.3);
  \draw[ultra thick,black] (.2,0) -- (.2,.3);
   \end{tikzpicture}   &$\ _3H^1_2(2,1)$ &$12$  \\
   \hline
    \rule{0ex}{3.5ex}
    \begin{tikzpicture}
    \draw[ultra thick, black] (0,0) rectangle (.25, .5);
   \end{tikzpicture}   &$\ _3P^2(3,2)$ &$8$  \\
   \hline
  \rule{0ex}{3.5ex}
  \begin{tikzpicture}
  \draw[ultra thick,black] (-.05,0) -- (-.05,.5);
   \draw[ultra thick,black] (.25,0) -- (.25,.5);
   \end{tikzpicture}   &$\ _3P^1(3,1)$ &$24$  \\
   \hline\hline
   \multicolumn{2}{ | l}
   {
   $\text{Sum}/48 = $ }& $96/48=2$ \\
   \hline
\end{tabular}

&\qquad
       \begin{tabular}{| c | c | c  |  }
       \hline
    ${\mathcal T}$ & Class Type & Summand  
      \\ \hline\hline
        \rule{0ex}{3.5ex}
        \begin{tikzpicture}
  \filldraw[black] (0,.4) circle (1.5pt);  
 \filldraw[black] (.2,.4) circle (1.5pt);
  \filldraw[black] (0,.2) circle (1.5pt);  
 \filldraw[black] (.2,.2) circle (1.5pt);
  \filldraw[black] (0,0) circle (1.5pt);  
 \filldraw[black] (.2,0) circle (1.5pt);
  \end{tikzpicture}
 & $Id$ & $800$     \\
 \hline
 \rule{0ex}{3.5ex}
  \begin{tikzpicture}
 \draw[ultra thick,black] (-.05,.5) -- (.25,.5);
 \draw[ultra thick,black] (-.05, .25) -- (.25,.25);
  \draw[ultra thick,black] (-.05,0) -- (.25,0);
   \end{tikzpicture}   &$\ _5C(1,2)$ &$160$  \\
   \hline
 \rule{0ex}{3.5ex}
 \begin{tikzpicture}
 \filldraw[black] (0,.5) circle (1.5pt);  
 \filldraw[black] (.2,.5) circle (1.5pt);
 \draw[ultra thick,black] (0,0) -- (0,.3);
  \draw[ultra thick,black] (.2,0) -- (.2,.3);
 \end{tikzpicture} & $\ _5H^1_2(2,1)$&$1200$  \\
        \hline
        \rule{0ex}{3.5ex}
         \begin{tikzpicture}
 \draw[ultra thick,black] (-.05,.5) -- (.25,.5);
 \draw[ultra thick,black] (0,0) -- (0,.3);
  \draw[ultra thick,black] (.2,0) -- (.2,.3);
   \end{tikzpicture}   &$\ _5H^1_2(2,1)$ &$240$  \\
   \hline
    \rule{0ex}{3.5ex}
    \begin{tikzpicture}
    \draw[ultra thick, black] (0,0) rectangle (.25, .5);
   \end{tikzpicture}    &$\ _5E(3,2)$ &$80$  \\
   \hline
  \rule{0ex}{3.5ex}
  \begin{tikzpicture}
  \draw[ultra thick,black] (-.05,0) -- (-.05,.5);
   \draw[ultra thick,black] (.25,0) -- (.25,.5);
   \end{tikzpicture}   &$\ _5E(3,1)$ &$400$  \\
   \hline\hline
   \multicolumn{2}{ | l}
   {
   $\text{Sum}/480 = $ }& $2880/480=6$ \\
   \hline
\end{tabular} \\
p=3 & \qquad p=5 
\end{array} 
$$

Computations for all admissible $R$-partitions, $R=3,4,5,6$, and $p=3,5,7$,   are given in Section~\ref{S:computations}.     
    Tables~\ref{Ta:fullcount} and \ref{Ta:elements} contain the explicit numerical data needed for  the computations.   
    
    \subsection{Surface automorphisms}\label{S:surfaces2} 
    Given $R$ and $p$, the sum $\sum_{\mathcal P} |O({\mathcal P})|$ over all admissible $R$-partitions can be interpreted as  the number of topological types of fully ramified ${\mathbb Z}_p^2$ actions in genus $g = 1 + Rp(p-1)/2 - p^2$, as described in Section~\ref{S:surfaces}.  Data of this type,  gleaned from the computations in Section~\ref{S:computations}, are   summarized in Table~\ref{Ta:topTypes} below.

\begin{table}[h]
\begin{center}
\begin{tabular}{| c || c | c || c |c  || c | c |}
  \hline
     & \multicolumn{2}{| c ||}{$p=3$} & \multicolumn{2}{|c ||}{$p=5$ }&\multicolumn{2}{| c |}{$p=7$}\\
     \hline
    ${\bf R}$ & Genus & $\#$ Top.Types & Genus & $\#$ Top. Types & Genus & $\#$ Top. Types \\
     \hline \hline
     {\bf $3$} & $1$&$1$ &$6$ & $1$ & $15$& $1$ \\
     \hline
     {\bf $4$} & $4$ & $2$ & $16$ & $4$ & $36$ & $6$ \\
     \hline
     {\bf $5$} & $7$ & $4$ & $26$ & $14$ & $57$ & $39$ \\
     \hline
{\bf $6$} & $10$ & $9$ & $36$ & $58$ & $78$ & $282$ \\
\hline
     \end{tabular}
     \caption{Topological types  of ${\mathbb Z}_p^2$ actions fully ramified over $R$ points, $p=3,5,7$.}\label{Ta:topTypes}
     \end{center}
     \end{table}

\section{Computations}\label{S:computations}
For each $R$-partition ${\mathcal P}$ and  prime $p$ listed below, the summands in brackets in \eqref{E:specialBurnside} are given with 
$|{\mathcal M}({\mathcal P})|$ listed first,  followed by triples consisting of a class type, a tiling ${\mathcal T}$, and the integer $t(\sigma)c(\sigma) |\text{Fix}_{\mathcal T}(\sigma)|$ for a representative element $\sigma$ in the class type.  The corresponding surface genus $g=1 + Rp(p-1)/2 - p^2$ is included for convenience.
\vskip .1 in
 \begin{description}
\item[$R=3$] \ 
 \begin{description}
 \item[$p=3, g=1$] {\bf 1 orbit}
 \begin{description}
 \item[$1^{[3]}$] $8$;\  $P^1(3,1)$ \begin{tikzpicture}
  \draw[ultra thick,black] (0,0) -- (0,.3);
   \end{tikzpicture}\  
   $16$;\  $H^1_2(2,1)$\ \begin{tikzpicture}
 \filldraw[black] (0,.34) circle (1.2pt);  
  \draw[ultra thick,black] (0,0) -- (0,.24);
   \end{tikzpicture}\   \ $24$; \ 
 $(8+16+24)/48 = {\bf 1}$.
 \end{description}
 \item[$p=5, g=6$] {\bf 1 orbit} 
  \begin{description}
 \item[$1^{[3]}$]  $80$;\  $E(3,1)$\  \begin{tikzpicture}
  \draw[ultra thick,black] (0,0) -- (0,.3);
   \end{tikzpicture}\  
   $160$;\  $H^1_2(2,1)$\ \begin{tikzpicture}
 \filldraw[black] (0,.34) circle (1.2pt);  
  \draw[ultra thick,black] (0,0) -- (0,.24);
   \end{tikzpicture}
   $240$; \ 
 $(80+160+240)/480 = {\bf 1}$.
 \end{description}

 \item[$p=7, g=15$] {\bf 1 orbit} 
  \begin{description}
 \item[$1^{[3]}$] $336$;\  $H^1_2(2,1)$\ \begin{tikzpicture}
 \filldraw[black] (0,.34) circle (1.2pt);  
  \draw[ultra thick,black] (0,0) -- (0,.24);
   \end{tikzpicture}
  \ $1008$;\  $H^3_3(3,1)$\ \begin{tikzpicture}
  \draw[ultra thick,black] (0,0) -- (0,.3);
   \end{tikzpicture}\  
   $672$;\  \ 
  $(336+1008+672)/ 2016 = {\bf 1}$.
 \end{description}

 \end{description}
 
 \item[$R=4$]\ 
 
 \begin{description}
 \item[$p=3, g=4$] {\bf 2 orbits} 
 \begin{description}
\item[$1^{[4]}$] 0;
 \item[$2^{[1]}1^{[2]}$] $24$;\  $H^1_2(2,1)$\ \begin{tikzpicture}
 \filldraw[black] (0,.34) circle (1.2pt);  
 \filldraw[black] (.13,.34) circle (1.2pt);
   \draw[ultra thick,black] (0,0) -- (0,.24);
   \end{tikzpicture}
  \ $24$; \ 
 $(24+24)/48 = {\bf 1}$.
 \item[$2^{[2]}$]  $6$;\  $C(1,2)$\ \begin{tikzpicture}
 \draw[ultra thick,black] (-.05, .16) -- (.16,.16);
  \draw[ultra thick,black] (-.05,0) -- (.16,0);
   \end{tikzpicture}\  
  \ $6$;\  $E(2,2)$\ \begin{tikzpicture}
   \draw[ultra thick,black] (0,0) rectangle (.2, .2);
  \end{tikzpicture}  
  $12$;\  $H^1_2(2,1)$\ \begin{tikzpicture}
 \draw[ultra thick,black] (0,0) -- (0,.26);
   \draw[ultra thick,black] (.13,0) -- (.13,.26);
 \end{tikzpicture} 
  \ $12$;\  $H^1_2(1,2)$\ \begin{tikzpicture}
 \filldraw[black] (0,.16) circle (1.2pt);  
 \filldraw[black] (.13,.16) circle (1.2pt);
  \draw[ultra thick,black] (-.05,0) -- (.17,0);
   \end{tikzpicture}\ 
  $12$; \\
  $(6+6 + 12+12+ 12)/48 = {\bf 1}$.

 \end{description}
 \item[$p=5, g=16$] {\bf 4 orbits}
 \begin{description}
 \item[$1^{[4]}$] $120$;\  $H^1_2(2,1)$\ \begin{tikzpicture}
  \draw[ultra thick,black] (0,.28) -- (0,.50);
\draw[ultra thick,black] (0,0) -- (0,.24);
   \end{tikzpicture}
  \ $120$;\  $H^2_4(4,1)$\ \begin{tikzpicture}
  \draw[ultra thick,black] (0,0) -- (0,.50);
   \end{tikzpicture}
  \ $240$; \ \ 
 $(120 + 120 + 240)/480 = {\bf 1}$.
  \item[$ 2^{[1]}1^{[2]}$]  $480$;\  $H^1_2(2,1)$\ \begin{tikzpicture}
 \filldraw[black] (0,.34) circle (1.2pt);  
 \filldraw[black] (.13,.34) circle (1.2pt);
   \draw[ultra thick,black] (0,0) -- (0,.24);
   \end{tikzpicture}
  \ $480$; \ \ 
   $(480 + 480)/480 = {\bf 2}$.
  \item[$2^{[2]}$]  $60$;\  $C(1,2)$\ \begin{tikzpicture}
 \draw[ultra thick,black] (-.05, .16) -- (.16,.16);
  \draw[ultra thick,black] (-.05,0) -- (.16,0);
   \end{tikzpicture}\  
  \ $60$; $H^1_2(2,1)$\ \begin{tikzpicture}
 \draw[ultra thick,black] (0,0) -- (0,.26);
   \draw[ultra thick,black] (.13,0) -- (.13,.26);
 \end{tikzpicture} 
  \ $120$;\  $H^1_2(2,1)$\ \begin{tikzpicture}
 \filldraw[black] (0,.16) circle (1.2pt);  
 \filldraw[black] (.13,.16) circle (1.2pt);
  \draw[ultra thick,black] (-.05,0) -- (.17,0);
    \end{tikzpicture}\ 
  $120$;\  $H^4_4(2,2)$ \begin{tikzpicture}
   \draw[ultra thick,black] (0,0) rectangle (.2, .2);
  \end{tikzpicture}  
  $120$; \\
  $(60 + 60 + 120 + 120 + 120)/480 = {\bf 1}$.
 \end{description}
 \item[$p=7, g=36$] {\bf 6 orbits} 
 \begin{description}
 \item[$1^{[4]}$] $1680$;\  $H^1_2(2,1)$\ \begin{tikzpicture}
  \draw[ultra thick,black] (0,.28) -- (0,.50);
\draw[ultra thick,black] (0,0) -- (0,.24);
   \end{tikzpicture}
  \  $1008$;\  $H^1_3(3,1)$\ \begin{tikzpicture}
 \filldraw[black] (0,.48) circle (1.2pt);  
  \draw[ultra thick,black] (0,0) -- (0,.34);
   \end{tikzpicture}
  \ $1344$; \ \ 
 $(1680 + 1008 + 1344)/2016 = {\bf 2}$.
   \item[$ 2^{[1]}1^{[2]}$]  $3024$; $H^1_2(2,1)$\ \begin{tikzpicture}
 \filldraw[black] (0,.34) circle (1.2pt);  
 \filldraw[black] (.13,.34) circle (1.2pt);
   \draw[ultra thick,black] (0,0) -- (0,.24);
   \end{tikzpicture}
  \ $3024$; \ \ 
   $(3024 + 3024)/2016 = {\bf 3}$.
  \item[$ 2^{[2]}$] $252$;  $C(1,2)$\ \begin{tikzpicture}
 \draw[ultra thick,black] (-.05, .16) -- (.16,.16);
  \draw[ultra thick,black] (-.05,0) -- (.16,0);
   \end{tikzpicture}\  
  \ $252$;\  $E(2,2)$ \begin{tikzpicture}
   \draw[ultra thick,black] (0,0) rectangle (.2, .2);
  \end{tikzpicture}  
  $504$;\  $H^1_2(2,1)$\ \begin{tikzpicture}
 \draw[ultra thick,black] (0,0) -- (0,.26);
   \draw[ultra thick,black] (.13,0) -- (.13,.26);
 \end{tikzpicture} 
  \ $504$;\   $H^1_2(2,1)$\ \begin{tikzpicture}
 \filldraw[black] (0,.16) circle (1.2pt);  
 \filldraw[black] (.13,.16) circle (1.2pt);
  \draw[ultra thick,black] (-.05,0) -- (.17,0);
    \end{tikzpicture}
  \ $504$; \\
 $(252+252+ 504+ 504+504)/2016 = {\bf 1}$.
 \end{description}

 \end{description}

 \item[$R=5$]\ 
 
 \begin{description}
 \item[$p=3, g=7$] {\bf 4 orbits}
 \begin{description}
 \item[$1^{[5]}$] $0$;
 \item[$2^{[1]}1^{[3]}$] $8$;\  $P^1(3,1)$\ \begin{tikzpicture}
 \filldraw[black] (0,.43) circle (1.2pt);  
  \filldraw[black] (.13,.43) circle (1.2pt); 
  \draw[ultra thick,black] (0,-.05) -- (0,.29);
   \end{tikzpicture}
  $16$;\ $H^1_2(2,1)$ \begin{tikzpicture}
  \filldraw[black] (0,.48) circle (1.2pt);
  \draw[ultra thick,black] (-.05,.34) -- (.16, .34);
 \draw[ultra thick,black] (0,0) -- (0,.24);
 \end{tikzpicture} 
 \ $24$; \ \ 
  $(8+16+24)/48 = {\bf 1}.$
 \item[$2^{[2]}1^{[1]}$] $24$;\  $H^1_2(2,1)$\ \begin{tikzpicture}
 \filldraw[black] (0,.34) circle (1.2pt);  
   \draw[ultra thick,black] (0,0) -- (0,.24);
   \draw[ultra thick,black] (.13,0) -- (.13,.24);
   \end{tikzpicture}
  \  \ $24$; \ \ 
 $(24+24)/48 = {\bf 1}.$
 \item[$3^{[1]}1^{[2]}$] $24$;\ $H^1_2(2,1)$ \ \begin{tikzpicture}
 \filldraw[black] (0,.34) circle (1.2pt);  
 \filldraw[black] (.13,.34) circle (1.2pt);
  \filldraw[black] (.26,.34) circle (1.2pt);
   \draw[ultra thick,black] (0,0) -- (0,.24);
   \end{tikzpicture}
  \ \ $24$; \ \ 
  $(24+24)/48 = {\bf 1}$.
 \item[$3^{[1]}2^{[1]}$] $24$;\ $H^1_2(2,1)$ \ \begin{tikzpicture}
 \filldraw[black] (0,.16) circle (1.2pt);  
 \filldraw[black] (.13,.16) circle (1.2pt);
 \filldraw[black] (.26,.16) circle (1.2pt);
  \draw[ultra thick,black] (-.05,0) -- (.17,0);
    \end{tikzpicture}\ 
  \ \ $24$; \ \ 
  $(24 + 24)/48 = {\bf 1}.$
  \end{description}
  
 \item[$p=5, g=26$] {\bf 14 orbits}
  \begin{description}
 \item[$1^{[5]}$] $264$;\ $P^1(5,1)$ \begin{tikzpicture}
 \draw[ultra thick,black] (0,0) -- (0,.59);
   \end{tikzpicture}
  \ $96$;\ $H^1_2(2,1)$ \begin{tikzpicture}
 \filldraw[black] (0,.59) circle (1.2pt);  
 \draw[ultra thick,black] (0,.28) -- (0,.50);
\draw[ultra thick,black] (0,0) -- (0,.24);
   \end{tikzpicture}
  $360$;\  $H^1_4(4,1)$\ \begin{tikzpicture}
 \filldraw[black] (0,.59) circle (1.2pt);  
 \draw[ultra thick,black] (0,0) -- (0,.50);
   \end{tikzpicture}
  \ $240$;\\
   $(264+ 96 + 360 + 240)/480 = {\bf 2}.$
 \item[$2^{[1]}1^{[3]}$] $1440$;\ $H^1_2(2,1)$ \ \begin{tikzpicture}
  \filldraw[black] (0,.48) circle (1.2pt);
  \draw[ultra thick,black] (-.05,.34) -- (.16, .34);
 \draw[ultra thick,black] (0,0) -- (0,.24);
 \end{tikzpicture} 
  \ $480$; \ \  
  $(1440+480)/480 = {\bf 4}.$
 \item[$2^{[2]}1^{[1]}$] $960$;\ $H^1_2(2,1)$\ \begin{tikzpicture}
 \filldraw[black] (0,.34) circle (1.2pt);  
   \draw[ultra thick,black] (0,0) -- (0,.24);
   \draw[ultra thick,black] (.13,0) -- (.13,.24);
   \end{tikzpicture}
  \  \ $480$; \ \ 
  $(960+480)/480 = {\bf 3}$.
 \item[$3^{[1]}1^{[2]}$] $960$;\ $H^1_2(2,1)$ \ \begin{tikzpicture}
 \filldraw[black] (0,.34) circle (1.2pt);  
 \filldraw[black] (.13,.34) circle (1.2pt);
  \filldraw[black] (.26,.34) circle (1.2pt);
   \draw[ultra thick,black] (0,0) -- (0,.24);
   \end{tikzpicture}
  \ \ $960$; \ \ 
  $(960+960)/480 = {\bf 4}.$
 \item[$3^{[1]}2^{[1]}$] $240$; \ $H^1_2(2,1)$ \ \begin{tikzpicture}
 \filldraw[black] (0,.16) circle (1.2pt);  
 \filldraw[black] (.13,.16) circle (1.2pt);
 \filldraw[black] (.26,.16) circle (1.2pt);
  \draw[ultra thick,black] (-.05,0) -- (.17,0);
    \end{tikzpicture}\ 
 \ \ $240$; \ \ 
  $(240 +240)/480 = {\bf 1}.$
  \end{description}
  
 \item[$p=7, g=57$] {\bf 39 orbits}
  \begin{description}
 \item[$1^{[5]}$] $9072$;\ $H^1_2(2,1)$ \begin{tikzpicture}
 \filldraw[black] (0,.59) circle (1.2pt);  
 \draw[ultra thick,black] (0,.28) -- (0,.50);
\draw[ultra thick,black] (0,0) -- (0,.24);
   \end{tikzpicture}
   $5040$;\ \ $(9072 + 5040)/2016 = {\bf 7}$.
 \item[$2^{[1]}1^{[3]}$] $25200$;\ $H^1_2(2,1)$ \ \begin{tikzpicture}
  \filldraw[black] (0,.48) circle (1.2pt);
  \draw[ultra thick,black] (-.05,.34) -- (.16, .34);
 \draw[ultra thick,black] (0,0) -- (0,.24);
 \end{tikzpicture} 
  \ $3024$;\  $H^1_3(3,1)$ \ \begin{tikzpicture}
 \filldraw[black] (0,.43) circle (1.2pt);  
  \filldraw[black] (.13,.43) circle (1.2pt); 
  \draw[ultra thick,black] (0,-.05) -- (0,.29);
   \end{tikzpicture}
  \ $4032$; \\
  $(25200 +3024+4032)/2016 = {\bf 16}.$
 \item[$2^{[2]}1^{[1]}$] $9072$;\ $H^1_2(2,1)$\ \begin{tikzpicture}
 \filldraw[black] (0,.34) circle (1.2pt);  
   \draw[ultra thick,black] (0,0) -- (0,.24);
   \draw[ultra thick,black] (.13,0) -- (.13,.24);
   \end{tikzpicture}
  \  \ $3024$;\ \ $(9072+3024)/2016={\bf 6}.$
 \item[$3^{[1]}1^{[2]}$] $8064$;\ $H^1_2(2,1)$ \ \begin{tikzpicture}
 \filldraw[black] (0,.34) circle (1.2pt);  
 \filldraw[black] (.13,.34) circle (1.2pt);
  \filldraw[black] (.26,.34) circle (1.2pt);
   \draw[ultra thick,black] (0,0) -- (0,.24);
   \end{tikzpicture}
  \ \ $8064$;\ \ $(8064 + 8064)/2016 = {\bf 8}.$
 \item[$3^{[1]}2^{[1]}$] $1344$;\ $H^1_2(2,1)$ \ \begin{tikzpicture}
 \filldraw[black] (0,.16) circle (1.2pt);  
 \filldraw[black] (.13,.16) circle (1.2pt);
 \filldraw[black] (.26,.16) circle (1.2pt);
  \draw[ultra thick,black] (-.05,0) -- (.17,0);
    \end{tikzpicture}\ 
  \ \ $1344$;\ $H^1_3(3, 1)$\ \begin{tikzpicture}
 \filldraw[black] (0,.16) circle (1.2pt);  
 \filldraw[black] (.13,.16) circle (1.2pt);
  \draw[ultra thick,black] (-.05,0) -- (.27,0);
    \end{tikzpicture}\ 
  \ $672$;\  $H^2_3(6,1)$ \ \begin{tikzpicture}
 \draw[ultra thick,black] (-.05, .16) -- (.16,.16);
  \draw[ultra thick,black] (-.05,0) -- (.27,0);
   \end{tikzpicture}\  
  \ $672$; \\
  $(1344+1344+672+672)/2016 = {\bf 2}.$
  \end{description}
  
 \end{description}
 
 \item[$R=6$]\ 
 
 \begin{description}
 \item[$p=3, g=10$] {\bf 9 orbits}
 \begin{description}
 \item[$1^{[6]}$] $0$; 
 \item[$2^{[1]}1^{[4]}$] $0$; 
 \item[$2^{[2]}1^{[2]}$] $24$;\ $H^1_2(2,1)$ \begin{tikzpicture}
 \filldraw[black] (0,.48) circle (1.2pt);  
 \filldraw[black] (.13,.48) circle (1.2pt);
 \draw[ultra thick,black] (-0.05,.34) -- (.16,.34);
 \draw[ultra thick,black] (0,0) -- (0,.26);
   \end{tikzpicture}\ 
   $24$;\ $(24 + 24)/48 = {\bf 1}$.
 \item[$3^{[1]}1^{[3]}$] $16$;\  $P^1(3,1)\ $\begin{tikzpicture}
 \filldraw[black] (0,.48) circle (1.2pt);  
 \filldraw[black] (.13,.48) circle (1.2pt);
 \filldraw[black] (.26,.48) circle (1.2pt);
 \draw[ultra thick,black] (0,0) -- (0,.34);
   \end{tikzpicture}\ 
   $32$;\ $(16 + 32)/48 = {\bf 1}$.
 \item[$4^{[1]}1^{[2]}$] $48$\ $H^1_2(2,1)\ $\begin{tikzpicture}
 \filldraw[black] (0,.34) circle (1.2pt);  
 \filldraw[black] (.13,.34) circle (1.2pt);
 \filldraw[black] (.26,.34) circle (1.2pt);
\filldraw[black] (.39,.34) circle (1.2pt);
  \draw[ultra thick,black] (0,0) -- (0,.26);
   \end{tikzpicture}\ 
     $48$;\ $(48+48)/48 = {\bf 2}$.
 \item[$3^{[1]}2^{[1]}1^{[1]}$] $48$;\ $48/48 = {\bf 1}$.
 \item[$2^{[3]}$] $12$;\ $C(1,2)$ \begin{tikzpicture}
 \draw[ultra thick,black] (-.05,.26) -- (.13,.26);
 \draw[ultra thick,black] (-.05, .13) -- (.13,.13);
  \draw[ultra thick,black] (-.05,0) -- (.13,0);
   \end{tikzpicture} 
   $4$;\ $P^1(3,1)$ \begin{tikzpicture}
  \draw[ultra thick,black] (0,0) -- (0,.3);
   \draw[ultra thick,black] (.13,0) -- (.13,.3);
   \end{tikzpicture} 
   $24$;\ $P^2(3,2)$ \begin{tikzpicture}
   \draw[ultra thick,black] (0,0) rectangle (.13, .3);
  \end{tikzpicture}  
   $8$;\ $H^1_2(2,1)$ \begin{tikzpicture}
 \filldraw[black] (0,.35) circle (1.2pt);  
 \filldraw[black] (.13,.35) circle (1.2pt);
 \draw[ultra thick,black] (0,0) -- (0,.26);
   \draw[ultra thick,black] (.13,0) -- (.13,.26);
 \end{tikzpicture} 
    $36$;\ $H^1_2(2,1)$ \begin{tikzpicture}
 \draw[ultra thick,black] (-0.05,.34) -- (.16,.34);
 \draw[ultra thick,black] (0,0) -- (0,.26);
   \draw[ultra thick,black] (.13,0) -- (.13,.26);
   \end{tikzpicture}
    $12$;\ \\
  $(12 + 4 + 24 + 8 + 36 + 12)/48 = {\bf 2}$.
 \item[$3^{[2]}$] $24$;\ $H^1_2(2,1)$  \begin{tikzpicture}
  \draw[ultra thick,black] (0,0) -- (0,.26);
   \draw[ultra thick,black] (.13,0) -- (.13,.26);
   \draw[ultra thick,black] (.26,0) -- (.26,.26);
 \end{tikzpicture} \ 
   $24$;\  
 $(24+24)/48 = {\bf 1}$.
 \item[$4^{[1]}2^{[1]}$] $12$;\ $C(1,2)$ \begin{tikzpicture}
 \draw[ultra thick,black] (.17,.16) -- (.34, .16);
 \draw[ultra thick,black] (-.05, .16) -- (.13,.16);
  \draw[ultra thick,black] (-.05,0) -- (.16,0);
   \end{tikzpicture}\  
    $12$;\ $H^1_2(2,1)$\ \begin{tikzpicture}
 \filldraw[black] (0,.16) circle (1.2pt);  
 \filldraw[black] (.13,.16) circle (1.2pt);
 \filldraw[black] (.26,.16) circle (1.2pt);
\filldraw[black] (.39,.16) circle (1.2pt);
  \draw[ultra thick,black] (-.05,0) -- (.15,0);
   \end{tikzpicture}\ 
     $12$;\  $H^1_2(2,1)$ \ \begin{tikzpicture}
 \filldraw[black] (0,.16) circle (1.2pt);  
 \filldraw[black] (.13,.16) circle (1.2pt);
  \draw[ultra thick,black] (-.05,0) -- (.15,0);
  \draw[ultra thick,black] (.18,0) -- (.35,0);
   \end{tikzpicture}\ 
    $12$;\  \\ 
  $(12+12+12+12)/48 = {\bf 1}$.
 \end{description}
 
 \item[$p=5, g=36$] {\bf 58 orbits}
  \begin{description}
 \item[$1^{[6]}$] $160$;\ $E(3,1)$ \begin{tikzpicture}
  \draw[ultra thick,black] (0,.28) -- (0,.50);
\draw[ultra thick,black] (0,0) -- (0,.24);
   \end{tikzpicture}
  \  $ 320$;\ $(160 + 320)/480 = {\bf 1}$.
 \item[$2^{[1]}1^{[4]}$] $3120$;\  $H^1_2(2,1)$\  
 \begin{tikzpicture}
 \filldraw[black] (0,.59) circle (1.2pt);  
 \filldraw[black] (.13,.59) circle (1.2pt);
 \draw[ultra thick,black] (0,.28) -- (0,.50);
\draw[ultra thick,black] (0,0) -- (0,.24);
   \end{tikzpicture}
  \ $960$;\ $H^1_2(2,1)$ \begin{tikzpicture}
 \draw[ultra thick,black] (-.05,.59) -- (.14,.59);
  \draw[ultra thick,black] (0,.28) -- (0,.50);
\draw[ultra thick,black] (0,0) -- (0,.24);
   \end{tikzpicture}
   \ $240$;\  $H^1_4(4,1)$\ \begin{tikzpicture}
 \filldraw[black] (0,.59) circle (1.2pt);  
 \filldraw[black] (.13,.59) circle (1.2pt);
 \draw[ultra thick,black] (0,0) -- (0,.50);
   \end{tikzpicture}
  \  $480$;\\
   $H^2_4(4,1)$  \begin{tikzpicture}
 \draw[ultra thick,black] (-.05,.57) -- (.14,.57);
  \draw[ultra thick,black] (0,0) -- (0,.50);
  \end{tikzpicture}
   \ $480$;\ 
  $(3120 + 960 + 240 + 480 + 480)/480 = {\bf 11}$.
 \item[$2^{[2]}1^{[2]}$] $5760$;\  \  $H^1_2(2,1)$ 
  \begin{tikzpicture}
 \filldraw[black] (0,.48) circle (1.2pt);  
 \filldraw[black] (.13,.48) circle (1.2pt);
 \draw[ultra thick,black] (-0.05,.34) -- (.16,.34);
 \draw[ultra thick,black] (0,0) -- (0,.26);
   \end{tikzpicture}\ 
\ \ $960$;\ $H^1_2(2,1)$ 
  \begin{tikzpicture}
  \draw[ultra thick,black] (0,.3) -- (0,.52);
 \draw[ultra thick,black] (.14,.3) -- (.14,.52);
\draw[ultra thick,black] (0,0) -- (0,.26);
   \end{tikzpicture}\ \ $480$;
   $(5760 + 960 +480)/480 = {\bf 15}.$
 \item[$3^{[1]}1^{[3]}$] $2880$;\ $2880/480 = {\bf 6}$.
 \item[$4^{[1]}1^{[2]}$] $1680$;\  $H^1_2(2,1)$
  \begin{tikzpicture}
 \filldraw[black] (0,.34) circle (1.2pt);  
 \filldraw[black] (.13,.34) circle (1.2pt);
 \filldraw[black] (.26,.34) circle (1.2pt);
\filldraw[black] (.39,.34) circle (1.2pt);
  \draw[ultra thick,black] (0,0) -- (0,.26);
   \end{tikzpicture}
  \ \  $1680$;\ $(1680 + 1680)/480 = {\bf 7}$.
 \item[$3^{[1]}2^{[1]}1^{[1]}$] $3840$; \  $3840/480 = {\bf 8}$.
  \item[$2^{[3]}$] $800$;\ $C(1,2)$ \begin{tikzpicture}
 \draw[ultra thick,black] (-.05,.26) -- (.13,.26);
 \draw[ultra thick,black] (-.05, .13) -- (.13,.13);
  \draw[ultra thick,black] (-.05,0) -- (.13,0);
   \end{tikzpicture} 
   $160$;\ $E(3,1)$ \begin{tikzpicture}
  \draw[ultra thick,black] (0,0) -- (0,.3);
   \draw[ultra thick,black] (.13,0) -- (.13,.3);
   \end{tikzpicture} 
   $400$;\ $E(3,2)$ \begin{tikzpicture}
   \draw[ultra thick,black] (0,0) rectangle (.13, .3);
  \end{tikzpicture}  
   $80$;\ $H^1_2(2,1)$ \begin{tikzpicture}
 \filldraw[black] (0,.35) circle (1.2pt);  
 \filldraw[black] (.13,.35) circle (1.2pt);
 \draw[ultra thick,black] (0,0) -- (0,.26);
   \draw[ultra thick,black] (.13,0) -- (.13,.26);
 \end{tikzpicture} 
    $1200$; \\
    $H^1_2(2,1)$ \begin{tikzpicture}
 \draw[ultra thick,black] (-0.05,.34) -- (.16,.34);
 \draw[ultra thick,black] (0,0) -- (0,.26);
   \draw[ultra thick,black] (.13,0) -- (.13,.26);
   \end{tikzpicture}
    $240$;\ 
  $(800 + 160 + 400 + 80 + 1200 + 240)/480 = {\bf 6}$.
 \item[$3^{[2]}$] $240$; $H^1_2(2,1)$\ \begin{tikzpicture}
  \draw[ultra thick,black] (0,0) -- (0,.26);
   \draw[ultra thick,black] (.13,0) -- (.13,.26);
   \draw[ultra thick,black] (.26,0) -- (.26,.26);
 \end{tikzpicture}
  \ $240$;\ $(240 + 240)/480 = {\bf 1}$. 
  \item[$4^{[1]}2^{[1]}$] $420$;\ $C(1,2)$ \begin{tikzpicture}
 \draw[ultra thick,black] (.17,.16) -- (.34, .16);
 \draw[ultra thick,black] (-.05, .16) -- (.13,.16);
  \draw[ultra thick,black] (-.05,0) -- (.16,0);
   \end{tikzpicture}\  
    $180$;\ $H^1_2(2,1)$\ \begin{tikzpicture}
 \filldraw[black] (0,.16) circle (1.2pt);  
 \filldraw[black] (.13,.16) circle (1.2pt);
 \filldraw[black] (.26,.16) circle (1.2pt);
\filldraw[black] (.39,.16) circle (1.2pt);
  \draw[ultra thick,black] (-.05,0) -- (.15,0);
   \end{tikzpicture}\ 
     $420$;\  $H^1_2(2,1)$ \ \begin{tikzpicture}
 \filldraw[black] (0,.16) circle (1.2pt);  
 \filldraw[black] (.13,.16) circle (1.2pt);
  \draw[ultra thick,black] (-.05,0) -- (.15,0);
  \draw[ultra thick,black] (.18,0) -- (.35,0);
   \end{tikzpicture}\ 
    $180$;\  \\ 
  $H^1_4(4,1)$ \ \begin{tikzpicture}
 \filldraw[black] (0,.16) circle (1.2pt);  
 \filldraw[black] (.13,.16) circle (1.2pt);
  \draw[ultra thick,black] (-.05,0) -- (.35,0);
   \end{tikzpicture}
  \  $120$;\ $H^2_4(4,1)$ \ \begin{tikzpicture}
   \draw[ultra thick,black] (-.05, .16) -- (.13,.16);
   \draw[ultra thick,black] (-.05,0) -- (.35,0);
   \end{tikzpicture}
  \  $120$;\ \\
  $(420 + 180 + 420 + 180 + 120 + 120)/480 = {\bf 3}$.
 \end{description}
 
 \item[$p=7, g=78$] {\bf 282 orbits}
  \begin{description}
 \item[$1^{[6]}$] $26544$; $H^1_2(2,1)$ \begin{tikzpicture}
 \draw[ultra thick,black] (0,.54) -- (0,.76);
  \draw[ultra thick,black] (0,.28) -- (0,.50);
\draw[ultra thick,black] (0,0) -- (0,.24);
   \end{tikzpicture}
  \  $1680$; $H^1_3(3,1)$  \begin{tikzpicture}
  \draw[ultra thick,black] (0,.38) -- (0,.76);
\draw[ultra thick,black] (0,0) -- (0,.34);
   \end{tikzpicture}
  \  $672$; $H^3_3(3,1)$  \begin{tikzpicture}
  \draw[ultra thick,black] (0,.38) -- (0,.76);
\draw[ultra thick,black] (0,0) -- (0,.34);
   \end{tikzpicture} \ $2016$;  $H^3_6(6,1)$\ \begin{tikzpicture}
\draw[ultra thick,black] (0,0) -- (0,.76);
   \end{tikzpicture}\ 
  $672$; $H^2_3(6,1)$ \begin{tikzpicture}
\draw[ultra thick,black] (0,0) -- (0,.76);
   \end{tikzpicture}\ 
  \ $672$; \ $(26544+1680+672 + 2016+672+672)/2016 = {\bf 16}.$
 \item[$2^{[1]}1^{[4]}$] $156240$; \ $H^1_2(2,1)$\ \begin{tikzpicture}
 \filldraw[black] (0,.59) circle (1.2pt);  
 \filldraw[black] (.13,.59) circle (1.2pt);
 \draw[ultra thick,black] (0,.28) -- (0,.50);
\draw[ultra thick,black] (0,0) -- (0,.24);
   \end{tikzpicture}
  \  \ $18144$;\ $H^1_2(2,1)$ \begin{tikzpicture}
 \draw[ultra thick,black] (-.05,.59) -- (.14,.59);
  \draw[ultra thick,black] (0,.28) -- (0,.50);
\draw[ultra thick,black] (0,0) -- (0,.24);
   \end{tikzpicture}
  \  $3024$; \\ 
  $(156240 +18144 + 3024)/2016 = {\bf 88}$. 
 \item[$2^{[2]}1^{[2]}$] $136080$; \  $H^1_2(2,1)$ 
  \begin{tikzpicture}
 \filldraw[black] (0,.48) circle (1.2pt);  
 \filldraw[black] (.13,.48) circle (1.2pt);
 \draw[ultra thick,black] (-0.05,.34) -- (.16,.34);
 \draw[ultra thick,black] (0,0) -- (0,.26);
   \end{tikzpicture}\ 
\ \ $9072$;\ $H^1_2(2,1)$ 
  \begin{tikzpicture}
  \draw[ultra thick,black] (0,.3) -- (0,.52);
 \draw[ultra thick,black] (.14,.3) -- (.14,.52);
\draw[ultra thick,black] (0,0) -- (0,.26);
   \end{tikzpicture}\ \ $6048$;\\
   $(136080 + 9072 + 6048)/2016 = {\bf 75}.$
 \item[$3^{[1]}1^{[3]}$] $67200$;\ $H^1_3(3,1)$ \begin{tikzpicture}
 \filldraw[black] (0,.48) circle (1.2pt);  
 \filldraw[black] (.13,.48) circle (1.2pt);
 \filldraw[black] (.26,.48) circle (1.2pt);
 \draw[ultra thick,black] (0,0) -- (0,.34);
   \end{tikzpicture}
  \ $10752$;\  $H^3_3(3,1)$ \begin{tikzpicture}
 \draw[thick,black] (-.05,.54) -- (.26,.54);
 \draw[ultra thick,black] (-.05,.44) -- (.26,.44);
 \draw[ultra thick,black] (0,0) -- (0,.34);
   \end{tikzpicture}
  \  $1344$;\ $H^3_3(3,1)$ \begin{tikzpicture}
 \draw[ultra thick,black] (-.05,.54) -- (.26,.54);
 \draw[thick,black] (-.05,.44) -- (.26,.44);
 \draw[ultra thick,black] (0,0) -- (0,.34);
   \end{tikzpicture}
  \  $1344$; \\
  $(67200 + 10752 + 1344 + 1344)/2016 = {\bf 40}.$ 
 \item[$4^{[1]}1^{[2]}$] $18144$;\  $H^1_2(2,1)$ \begin{tikzpicture}
 \filldraw[black] (0,.34) circle (1.2pt);  
 \filldraw[black] (.13,.34) circle (1.2pt);
 \filldraw[black] (.26,.34) circle (1.2pt);
\filldraw[black] (.39,.34) circle (1.2pt);
  \draw[ultra thick,black] (0,0) -- (0,.26);
   \end{tikzpicture}
  \ \  $18144$;\ $(18144 + 18144)/2016 = {\bf 18}.$

 \item[$3^{[1]}2^{[1]}1^{[1]}$] $48384$; \  $(48384)/2016 = {\bf 24}.$
 
 \item[$2^{[3]}$] $10584$;\  $C(1,2)$\ \begin{tikzpicture}
 \draw[ultra thick,black] (-.05,.26) -- (.13,.26);
 \draw[ultra thick,black] (-.05, .13) -- (.13,.13);
  \draw[ultra thick,black] (-.05,0) -- (.13,0);
   \end{tikzpicture} 
  \  $1512$;\ $H^1_2(2,1)$\  \begin{tikzpicture}
 \filldraw[black] (0,.35) circle (1.2pt);  
 \filldraw[black] (.13,.35) circle (1.2pt);
 \draw[ultra thick,black] (0,0) -- (0,.26);
   \draw[ultra thick,black] (.13,0) -- (.13,.26);
 \end{tikzpicture} 
  \  $10584$;\ $H^1_2(2,1)$\ \begin{tikzpicture}
 \draw[ultra thick,black] (-0.05,.34) -- (.16,.34);
 \draw[ultra thick,black] (0,0) -- (0,.26);
   \draw[ultra thick,black] (.13,0) -- (.13,.26);
   \end{tikzpicture}
 \   $1512$;\ $H^1_3(3,1)$ \ \begin{tikzpicture}
  \draw[ultra thick,black] (0,0) -- (0,.3);
   \draw[ultra thick,black] (.13,0) -- (.13,.3);
   \end{tikzpicture} 
  \ $672$;\ $H^3_3(3,1)$ \ \begin{tikzpicture}
  \draw[ultra thick,black] (0,0) -- (0,.3);
   \draw[ultra thick,black] (.13,0) -- (.13,.3);
   \end{tikzpicture} 
  \  $2352$;\ $H^2_6(3,2)$\ \begin{tikzpicture}
   \draw[ultra thick,black] (0,0) rectangle (.13, .3);
  \end{tikzpicture}  
 \   $672$;\ $H^6_6(3,2))$ \ \begin{tikzpicture}
   \draw[ultra thick,black] (0,0) rectangle (.13, .3);
  \end{tikzpicture}  
 \  $336$; \\
$(10584+1512+10584+1512+672+2352+672+ 336)/2016 = {\bf 14}.$ 
 
  \item[$3^{[2]}$] $1792$;\ $C(1,3)\ $\begin{tikzpicture}
 \draw[ultra thick,black] (-.05, .16) -- (.23,.16);
  \draw[ultra thick,black] (-.05,0) -- (.23,0);
   \end{tikzpicture}
  \  $224$;\ $H^1_2(2,1)$\ \begin{tikzpicture}
  \draw[ultra thick,black] (0,0) -- (0,.26);
   \draw[ultra thick,black] (.13,0) -- (.13,.26);
   \draw[ultra thick,black] (.26,0) -- (.26,.26);
 \end{tikzpicture}
  \  $1344$;\ $H^1_3(3,1)$ \ \begin{tikzpicture}
 \filldraw[black] (0,.16) circle (1.2pt);  
 \filldraw[black] (.13,.16) circle (1.2pt);
  \filldraw[black] (.26,.16) circle (1.2pt);
  \draw[ultra thick,black] (-.05,0) -- (.35,0);
   \end{tikzpicture}
  \  $1792$; \\
   $H^3_3(3,1)$ \begin{tikzpicture}
 \draw[ultra thick,black] (-.05, .16) -- (.23,.16);
  \draw[ultra thick,black] (-.05,0) -- (.23,0);
   \end{tikzpicture}
  \  $224$;\ $H^3_6(2,3)$\ \begin{tikzpicture}
   \draw[ultra thick,black] (0,0) rectangle (.3, .16);
  \end{tikzpicture}  
  \   $672$;\\ 
  $(1792+224+1344+1792+ 224 + 672)/2016 = {\bf 3}.$
 
  \item[$4^{[1]}2^{[1]}$] $3024$; \ $C(1,2)$\ \begin{tikzpicture}
 \draw[ultra thick,black] (.17,.16) -- (.34, .16);
 \draw[ultra thick,black] (-.05, .16) -- (.13,.16);
  \draw[ultra thick,black] (-.05,0) -- (.16,0);
   \end{tikzpicture}
  \  $1008$;\ $H^1_2(2,1)$ \begin{tikzpicture}
 \filldraw[black] (0,.16) circle (1.2pt);  
 \filldraw[black] (.13,.16) circle (1.2pt);
 \filldraw[black] (.26,.16) circle (1.2pt);
\filldraw[black] (.39,.16) circle (1.2pt);
  \draw[ultra thick,black] (-.05,0) -- (.15,0);
   \end{tikzpicture}
  \ \  $ 3024$;\ $H^1_2(2,1)$ \ \begin{tikzpicture}
 \filldraw[black] (0,.16) circle (1.2pt);  
 \filldraw[black] (.13,.16) circle (1.2pt);
  \draw[ultra thick,black] (-.05,0) -- (.15,0);
  \draw[ultra thick,black] (.18,0) -- (.35,0);
   \end{tikzpicture}
  \    $1008$; \\
  $(3024+1008+3024+1008)/2016 ={\bf 4}.$

 \end{description}
 \end{description}
 \end{description}

\section*{Appendix:  Proof of Theorem~\ref{T:multiPart}}\label{S:appendix}
  The  {\it $q$-binomial coefficients} are 
 \begin{equation}\label{Gaussnum}
 \begin{bmatrix}
 n \\
 k
 \end{bmatrix}_q= \dfrac{(q^{n}-1)(q^{n-1}-1)\cdot \dots \cdot (q^{n-k+1}-1)}{(q^k-1)(q^{k-1}-1)\cdot \dots \cdot (q-1)} 
  \end{equation}
 where $n \geq k$ are positive integers.    If $q$ is a prime power, $[\begin{smallmatrix} n \\ k \end{smallmatrix}]_q$ is the number of $k$-dimensional subspaces of an $n$-dimensional vector space over the field with $q$ elements.     It is easily verified that 
 \begin{equation}\label{E:GaussIdentity}
  \begin{bmatrix}
 n \\
 k
 \end{bmatrix}_q=  \begin{bmatrix}
 n \\
 n-k
 \end{bmatrix}_q.
 \end{equation}
 Though defined as a rational functions, the $q$-binomial coefficients turn out to be polynomials  with integer coefficients,
 $$
  \begin{bmatrix}
 n \\
 k
 \end{bmatrix}_q = \sum_{l=0}^{k(n-k)}t_lq^l,
 $$
 where $t_l$ is equal to the number of partitions of $1 \leq l \leq k(n-k)$ into at most $k$ parts of size at most $n-k$.  A proof of this statement, and further details on the $q$ binomial coefficients,  can be found in, e.g.,  \cite{vLW}, Ch. 24.

 In accordance with the remarks preceding the statement of the theorem,  we may assume the weights $\omega_1 \equiv \omega_2 \equiv \dots \equiv \omega_m \equiv 1$ mod $p$.

 Let $a_{1j}\geq 0$ be the number of parts equal to $j$, $1 \leq j \leq p-1$,  in the first partition $r_1 + r_2 + \dots + r_{\mu_1}$ in \eqref{E:pseudoPart}.  Similarly let $a_{2j}$, $\dots$, $a_{mj}$ be the numbers of parts equal to $j$ in the remaining $m-1$ partitions comprising \eqref{E:pseudoPart}.  Then the  $m \times (p-1)$ matrix $[a_{ij}]$ satisfies the conditions 
     \begin{equation}\label{E:matCond1}
      \sum_{j=1}^{p-1} a_{ij}= \mu_i,\quad i = 1, \dots, m,
      \end{equation}
            and 
     \begin{equation}\label{E:matCond2}
     \sum_{i=1}^m  \sum_{j=1}^{p-1} ja_{ij} \equiv \alpha \pmod p,
     \end{equation}
     and specifies the multi-partition \eqref{E:pseudoPart} uniquely.  Thus determination of $n_{\alpha}$ is equvalent to  counting $m \times (p-1)$ matrices $[a_{ij}]$ of nonnegative integers satisfying  \eqref{E:matCond1} and \eqref{E:matCond2}.
     
     First suppose the matrices have an extra initial column $a_{i0}$, $i = 1, \dots, m$ and that  
     the expanded matrices satisfy the conditions analogous to \eqref{E:matCond1} and \eqref{E:matCond2}, namely,
     \begin{equation}\label{E:fullsum}
      \sum_{j=0}^{p-1} a_{ij}= \mu_i,\quad   i = 1, \dots, m,
      \end{equation}
      and 
            \begin{equation}\label{E:fullweight}
      \sum_{i=1}^m  \sum_{j=0}^{p-1}ja_{ij}  \equiv \alpha \pmod p,
      \end{equation}
      the latter being, in fact, identical to \eqref{E:matCond2}.  
 Let $\tilde n_{ \alpha}$ be the number of  expanded matrices satisfying \eqref{E:fullsum} and \eqref{E:fullweight}.      As a first step, we show that 
    \begin{equation}\label{E:0claim}
  (\tilde n_0, \tilde n_1, \tilde n_2, \dots, \tilde n_{p-1}) = \begin{cases}\biggl(\frac Ep, \frac Ep, \dots, \frac Ep \biggr) &\text{ if \  $\exists i,  \mu_i \not\equiv 0$ mod $p$;} \\
  \biggl(\frac{E-1}{p} + 1, \frac{E-1}{p}, \frac{E-1}{p}, \dots, \frac{E-1}{p}\biggr) & \text{otherwise,}\\
 \end{cases}
 \end{equation} 
 where 
$E=\prod_{i=1}^m e_{\mu_i}$, and $e_{\mu_i}= \binom{\mu_i+p-1}{\mu_i}$.  
$e_{\mu_i}$ is the number of ways of populating the $i$th row of the expanded matrix.  ($\mu_i$ `balls' distribute themselves over $p
 $ `boxes'  labelled $0, 1, \dots, p-1$.)     Hence the total number of matrices is $\prod_i e_{\mu_i} = E$.  Let $Y_{ \alpha}$ be the subset of matrices with weighted sum $\sum_{i=1}^m \sum_{j=0}^{p-1} ja_j$  in the specified residue class $\alpha$.   To prove the first case of \eqref{E:0claim}, apply the right shift map $\sigma$ to a row, say, the $i$th row, having $\mu_i \not\equiv 0$ mod $p$.
 \begin{gather}\notag
 \sigma: a_{ij} \mapsto a_{i\ j+1}, \quad 
 j=0, \dots p-2,  \\ \sigma: a_{i\ p-1} \mapsto a_{i0}.
 \end{gather}
Extend it by the identity map on the other rows. Then
 $$\sum_{i=1}^m \sum_{j=0}^{p-1}  j \sigma(a_{ij}) \equiv \sum_{i=1}^m \sum_{j=0}^{p-1} j a_{ij} +  \mu_i \pmod p.$$ 
 Construct the inverse map in the obvious way using a left shift map.  
We obtain a bijection between $Y_{\alpha}$ and $Y_{ \alpha +\mu_i}$.     Composing this bijection with itself $k$ times gives a bijection between $Y_{\alpha}$ and $Y_{\alpha +k\mu_i}$.   The numbers $\alpha + k\mu_i $, $k=0, 1, \dots, p-1$ determine the  complete set of residue classes mod $p$; this yields the first case of \eqref{E:0claim}.   To prove the second case of \eqref{E:0claim}, apply the following map, for  $\beta \not\equiv 0$ mod $p$, 
\begin{equation}\label{E:multiplier}
\phi_{\beta}: (a_{i0}, a_{i1}, \dots, a_{i\ p-1}) \mapsto (a_{i(\beta 0)}, a_{i(\beta 1)}, a_{i(\beta 2)}, \dots, a_{i(\beta (p-1))})
\end{equation}
to all rows of the matrix.   This 
 defines a bijection $\phi_\beta: Y_{\alpha} \rightarrow Y_{\beta\alpha}$ with inverse $\phi_{\beta^{-1}}: Y_{\beta\alpha} \rightarrow Y_{\alpha}$, as can be seen from
 $$
 \sum_{i=1}^m  \sum_{j=0}^{p-1} \beta j a_{ij} = \beta \sum_{i=1}^m \sum_{j=0}^{p-1} j a_{ij}.
 $$
  From this we conclude again that  $\tilde n_{ 1}=\tilde n_{ \alpha}$ for all $\alpha \not\equiv 0$ mod $p$, now including the case in which all $\mu_i \equiv 0$ mod $p$.   The second case of \eqref{E:0claim} will follow if we show that $\tilde n_{0} = \frac{E-1}{p}+1$.

The presence of an initial column  ($a_{i0}, i=1, \dots m$), with weight $j=0$, allows the $i$th row of $[a_{ij}]$ to be interpreted as a partition of $0 \leq l=\sum_{j=0}^{p-1}j a_{ij} \leq \mu_i(p-1)$ into at most $\mu_i$ positive parts of size at most $p-1$.   Hence the number of ways of populating the $i$th row with  weighted sum $l$ is the coefficient of $q^l$ in the polynomial form of the $q$-binomial coefficient
$$
\begin{bmatrix}
 \mu_i + p-1 \\
 \mu_i
 \end{bmatrix}_q = \sum_{l=0}^{\mu_i (p-1)} t_{l ,i}q^l.
$$
 It follows that the number of 
matrices $[a_{ij}]$ with total weighted sum
$$
L=\sum_{i=1}^m\sum_{j=0}^{p-1} ja_{ij}, \quad 0 \leq L \leq R(p-1), \quad R=\sum_i\mu_i,
$$
is the coefficient of $q^L$ in the polynomial product
$$
\prod_{i=1}^m \begin{bmatrix}
 \mu_i + p-1 \\
 \mu_i
 \end{bmatrix}_q = \sum_{L=0}^{R(p-1)} T_L q^L,
$$
where
$$
T_L= \sum_{l_1 + l_2 + \dots +  l_m =L}t_{l_i, i}.
$$
The number we seek is 
$$
\tilde n_{{\mathcal P}, 0}=\sum_{p \mid L} T_L.
$$
To pick out the desired coefficients, we use \begin{equation}\notag
 \sum_{j=0}^{p-1} (e^{2\pi i/p})^{jm} = \begin{cases} p & \text{if} \quad p\mid m \\ 0 & \text{if}\quad p \nmid m
\end{cases},
\end{equation}
where $e$ is the base of natural logarithms and $i$ is the imaginary unit (see \cite{Ap}, Theorem~8.1).  
It follows that 
\begin{equation}\notag p \tilde n_{0}=p\sum_{p\mid L}T_L  = \sum_{j=0}^{p-1}\sum_{L=0}^{R(p-1)} T_L  (e^{2\pi i/p})^{jL}.
\end{equation}
The first summand in the outer summation ($j=0$) is 
$$\sum_{L=0}^{R(p-1)} T_L = E,$$
and the remaining summands ($j=1, 2, \dots, p-1$) are all reorderings of the summation
$$\sum_{L=0}^{R(p-1)} T_L (e^{2\pi i/p})^{L}= \prod_{i=1}^m \begin{bmatrix}
 \mu_i + p-1 \\
 \mu_i
 \end{bmatrix}_{e^{2\pi i/p}} = 1^m = 1,$$
 using \eqref{Gaussnum}, \eqref{E:GaussIdentity},  and the assumption that each $\mu_i \equiv 0 \pmod p$. 
Thus 
\begin{equation}\label{E:0class}\notag
p\tilde n_{0}=E+p-1,
\end{equation}
 which completes the proof of second case of \eqref{E:0claim}.   
 
 We now return to the conditions in the statement of the theorem,  by assuming the extra initial column of the expanded matrix $[a_{ij}]$ consists entirely of $0$'s, i.e., $a_{i0}=0, i = 1, \dots, m$.  
The number of ways of choosing the $i$th row of such a matrix is 
$$b_{\mu_i} = e_{\mu_i} - e_{\mu_i-1} = \binom{\mu_i +p-2}{\mu_i}.$$
If both $\mu_i$ and $\mu_i - 1 \not\equiv 0$ mod $p$ (i.e., $\mu_i \not\equiv 0,1$ mod $p$), \eqref{E:0claim}  
 applied to the singletons ${\mathcal P} = (\mu_i)$ and ${\mathcal P} = (\mu_i-1)$, implies that there are 
 $$\frac{e_{\mu_i}-e_{\mu_i-1}}{p} = \frac{b_{\mu_i}}{p}$$
 ways of populating the $i$th row $[0, a_{i1} \dots, a_{ij} \dots a_{i\ p-1}]$ such that the weighted sum $\sum_{j=1}^{p-1} j a_{ij}$ belongs to any given residue class $\alpha$ mod $p$.  Then as long as ${\mathcal P}$ contains an element $\mu_i \not\equiv 0,1$ mod $p$, 
  \begin{equation}\notag
(n_{{\mathcal P}, 0}, n_{{\mathcal P}, 1}, \dots, n_{{\mathcal P}, p-1}) = \biggl(\frac{B_{\mathcal P}}{p}, \frac{B_{\mathcal P}}{p}, \dots, \frac{B_{\mathcal P}}{p}\biggr).
\end{equation}
 Even without this restriction on ${\mathcal P}$,   we have $n_{{\mathcal P},  \alpha}=n_{{\mathcal P},  \beta}$ if $\alpha, \beta \ne 0$, because the map \eqref{E:multiplier} fixes the initial column (of $0$'s).    Thus in the remaining case,  in which ${\mathcal P}$ consists entirely of elements $\mu_i \equiv 0$ or $\mu_i\equiv 1$ mod $p$, there exist nonnegative integers $\Omega_{{\mathcal P}}$ and $\zeta_{{\mathcal P}}$, satisfying $\Omega_{\mathcal P} +(p-1)\zeta_{\mathcal P}=B_{\mathcal P}$,  such that
$$(n_{{\mathcal P}, 0}, n_{{\mathcal P}, 1}, \dots, n_{{\mathcal P}, p-1})  = (\Omega_{\mathcal P}, \zeta_{\mathcal P}, \dots, \zeta_{\mathcal P}).$$ 
 To obtain expressions for $\Omega_{\mathcal P}$ and $\zeta_{\mathcal P}$,  relabel 
\begin{equation}\label{E:all01}
{\mathcal P}= (\mu_1, \mu_2, \dots, \mu_s, \nu_1, \nu_2, \dots, \nu_{t}), \quad s+t=m, \quad s,t \geq 0,
\end{equation}
so that  $\mu_i \equiv 0$ mod $p$ and  $\nu_j\equiv 1$ mod $p$.  Assume all the $\mu_i$ are nonzero.     
Consider the $m$-tuples $\tilde{\mathcal P}$ that can be formed from \eqref{E:all01} by reducing some, none,  or all elements by 1.  All but one of these tuples contain an element $\not\equiv 0$ mod $p$.  The exception is  
\begin{equation}\notag
{\mathcal P}_0 = (\mu_1, \mu_2, \dots, \mu_s, \nu_1-1, \dots, \nu_t-1).
\end{equation} 
  Applying  \eqref{E:0claim} to each $\tilde {\mathcal P} \neq {\mathcal P}_0$, we obtain 
  $$(n_{{\mathcal P}, 0}, n_{{\mathcal P}, 1}, \dots, n_{{\mathcal P}, p-1})\quad = \quad \biggl(\dfrac{\tilde E}{p}, \dfrac{\tilde E}{p}, \dots, \dfrac{\tilde E}{p} \biggr),$$
  for some $\tilde E$. 
 On the other hand,  applying \eqref{E:0claim} to $\tilde {\mathcal P}_0$ we obtain
 \begin{equation}\label{E:C'dist}
 (n_{{\mathcal P}, 0}, n_{{\mathcal P}, 1}, \dots, n_{{\mathcal P}, p-1})\quad = \quad \biggl(\frac{E_0-1}{p} + 1, \frac{E_0-1}{p}, \frac{E_0-1}{p}, \dots, \frac{E_0-1}{p}\biggr),
 \end{equation}
 where  
 $E_0= e_{\mu_1}e_{\mu_2} \dots e_{\mu_t} e_{\nu_1-1} e_{\nu_2-1} \dots e_{\nu_s - 1}$.  The numbers $\tilde E$, $E_0$ are (up to sign) the summands in the expansion of 
$$B_{\mathcal P}=\prod_{i=1}^s (e_{\mu_i}-e_{\mu_i-1})\prod_{j=1}^t(e_{\nu_j}-e_{\nu_j-1}). $$
  $E_0$ appears with sign $(-1)^t$, so $B_{\mathcal P}+(-1)^{t+1}E_0$ is the sum of all the other (signed) terms, each of which counts a set of matrices whose weighted sums, by  \eqref{E:0claim},  are equidistributed mod $p$.  $E_0$ counts a set of matrices whose weighted sums are distributed according to \eqref{E:C'dist}.  From this and the fact that $b_{\mu_i} = e_{\mu_i} - e_{\mu_i-1},$  we deduce that
  $$\zeta_{\mathcal P}=n_{{\mathcal P}, 1} =n _{{\mathcal P}, 2}= \dots = n_{{\mathcal P}, p-1} = \dfrac{B_{\mathcal P} + (-1)^{t+1}E_0}{p}+(-1)^t\dfrac{E_0-1}{p} $$
  while
  $$\Omega_{\mathcal P}=n_{{\mathcal P}, 0} = \dfrac{B_{\mathcal P} + (-1)^{t+1}E_0}{p}+(-1)^t\biggl(\dfrac{E_0-1}{p} +1\biggr). $$
  Terms involving $E_0$ drop out of both expressions, and we find  $\Omega_{\mathcal P} = W_{\mathcal P}$ and $\zeta_{\mathcal P} = Z_{\mathcal P}$ as in the statement of the theorem.

\subsection*{Acknowledgements}

Thanks to the referee of an earlier version of this paper, who  pointed out a significant omission.    Thanks to Cormac O'Sullivan for help with the proof of Theorem~\ref{T:multiPart}.

\end{document}